\documentclass[11pt]{amsart}
\usepackage{amsmath,amssymb,url}
\usepackage{tikz}
\usepackage{enumitem}

\newtheorem{theorem}{Theorem}

\newtheorem{corollary}[theorem]{Corollary}

\theoremstyle{definition}

\newtheorem{remark}[theorem]{Remark}

\newcommand{\NN}{\mathbb{N}} 
\newcommand{\ZZ}{\mathbb{Z}} 
 
\newcommand{\RR}{\mathbb{R}}

\newcommand{\A}{\mathrm{A}}

\newcommand{\E}{\mathrm{E}}

\newcommand{\V}{\mathrm{V}}

\renewcommand{\H}{\mathrm{H}}

\newcommand{\Gold}{{\mathrm{G}}}
\newcommand{\DjM}{{\mathrm{DjM}}}

\newcommand{\CT}{\mathrm{CT}}
\newcommand{\Aut}{\mathrm{Aut}}
\newcommand{\Alt}{\mathrm{Alt}}
\newcommand{\Sym}{\mathrm{Sym}}
\newcommand{\GL}{\mathrm{GL}}
\newcommand{\GF}{{\mathrm{GF}}}
\newcommand{\PSL}{{{\mathrm{PSL}}}}  
\newcommand{\PGL}{{{\mathrm{PGL}}}}
\newcommand{\PSU}{{{\mathrm{PSU}}}}
\newcommand{\PG}{{{\mathrm{PG}}}}
\newcommand{\rK}{{{\mathrm{K}}}}
\newcommand{\CSS}{{{\mathrm{CSS}}}}
\newcommand{\CAT}{{{\mathrm{CAT}}}}
\newcommand{\Cos}{\mathrm{Cos}}
\newcommand{\BiCos}{\mathrm{BiCos}}
\newcommand{\TypeSemicolon}{\mathrm{Type}\colon\> }

\newcommand{\Tree}{{\mathcal{T}}}
\newcommand{\cH}{{\mathcal{H}}}
\newcommand{\cG}{{\mathcal{G}}}
\newcommand{\cN}{{\mathcal{N}}}
\newcommand{\cP}{{\mathcal{P}}}

\newcommand{\cC}{{\mathcal{C}}}

\newcommand{\Ga}{{\Gamma}}

\newcommand{\tGa}{{\tilde{\Gamma}}}
\newcommand{\tG}{{\tilde{G}}}

\newcommand{\tK}{{\tilde{K}}}
\newcommand{\tA}{{\tilde{A}}}
\newcommand{\tg}{{\tilde{g}}}
\newcommand{\tv}{\tilde{v}}
\newcommand{\te}{\tilde{e}}

\newcommand{\norml}{\trianglelefteq}

\newcommand{\cal}{\mathcal}

\begin{document}

\title{Edge-transitive cubic graphs: Cataloguing and Enumeration}

\author[Marston Conder]{Marston Conder}
\address{Marston Conder,\newline
Department of Mathematics, University of Auckland,\newline
 38 Princes Street, Auckland 1010, New Zealand}
 \email{m.conder@auckland.ac.nz}
 
	\author[P.~Poto\v{c}nik]{Primo\v{z} Poto\v{c}nik}
	\address{Primo\v{z} Poto\v{c}nik,\newline
	Faculty of Mathematics and Physics,
	University of Ljubljana, \newline
	Jadranska ulica 19, 1000 Ljubljana, Slovenia; \newline
	also affiliated with\newline
	Institute of Mathematics, Physics, and Mechanics, \newline
	Jadranska ulica 19, 1000 Ljubljana, Slovenia}
	\email{primoz.potocnik@fmf.uni-lj.si}



\begin{abstract}
This paper deals with finite cubic ($3$-regular) graphs whose automorphism group acts transitively on the 
edges of the graph. Such graphs split into two broad classes, namely arc-transitive and semisymmetric 
cubic graphs, and then these divide respectively into $7$ types (according to a classification 
by Djokovi\'c and Miller (1980)) and $15$ types (according to a classification by Goldschmidt(1980)), 
in terms of certain group amalgams. 
Such graphs of small order were previously known up to orders $2048$ and $768$, respectively, 
and we have extended each of the two lists of all such graphs up to order $10000$. 
Before describing how we did that, we carry out an analysis of the $22$ amalgams, to show which 
of the finitely-presented groups associated with the $15$ Goldschmidt amalgams can be faithfully 
embedded in one or more of the other $21$ (as subgroups of finite index), complementing what is 
already known about such embeddings of the $7$ Djokovi\'c-Miller groups in each other.  
We also give an example of a graph of each of the $22$ types, and in most cases, describe the smallest 
such graph, and we then use regular coverings to prove that there are infinitely many examples of each type.  
Finally, we discuss the asymptotic enumeration of the graph orders, proving that if $f_\cC(n)$ is the 
number of cubic edge-transitive graphs of type $\cC$ on at most $n$ vertices, then there exist
positive real constants $a$ and $b$ and a positive integer $n_0$ such that   
$n^{a \log(n)} \le f_\cC(n) \le    n^{b \log(n)}$ for all $\ n\ge n_0$.
\end{abstract}

\maketitle

\section{Introduction}
\label{sec:Intro}

This paper deals with {\em cubic edge-transitive graphs}, or more precisely, $3$-regular graphs whose automorphism group acts transitively on the edges of the graph.
The reader may assume throughout the paper that all the graphs discussed are simple and connected, 
and moreover, with the exception of the (infinite) $3$-valent tree $\Tree_3$, all of them are finite.

In general, every edge-transitive regular graph with odd valency has one of two symmetry types, 
depending on whether its automorphism group acts transitively on the vertices or not.
If the action of the automorphism group on the vertices is transitive, then it is also transitive on the arcs 
(where an {\em arc} is defined as  an ordered pair of adjacent vertices), and in that case the graph is 
said to be {\em arc-transitive}, or sometimes {\em symmetric}.
On the other hand, if the action is intransitive on the vertices, then the graph is said to be {\em semisymmetric}, 
and in that case the graph is bipartite with the two vertex-orbits forming the bipartition.

The study of cubic edge-transitive graphs is one the oldest topics in algebraic graph-theory,
going back to Tutte's ingenious discovery and proof of an upper bound of $32$ on the order of 
a finite edge-stabiliser in the automorphism group of an arc-transitive cubic graph \cite{T47},
and an analogous theorem of Goldschmidt \cite{Gold80} for semisymmetric cubic graphs.
%
Tutte's and Goldschmidt's work had a profound impact not only on this particular branch of graph theory,
 but also on the development of local analysis and amalgams in group theory. 
The latter also played (and still do play) an important role in the Classification of Finite Simple Groups project.
 
Research on cubic edge-transitive graphs has in large part been driven and facilitated by the amazing and famous 
Foster census \cite{Foster}, which contained information on almost all of the cubic 
arc-transitive graphs on up to $512$ vertices, and its later completion and extension up to $768$ vertices, 
by Conder and Dobcs{\' a}nyi \cite{CD}, later complemented by the authors of the current paper jointly with Malni{\v c} and Maru{\v s}i{\v c} \cite{CMMP} in producing a list of all semisymmetric cubic graphs 
on up to $768$ vertices. 
The list in \cite{CD} was extended further in 2006 by the first author of the current paper, up to order $2048$, 
and later in 2011 up to order $10000$ (see \cite{ATConder}). 
These lists contained several previously undiscovered small graphs. 
The third smallest graph on the list in \cite{CMMP}, on 112 vertices, was the subject of a separate paper \cite{CMMPP}. 

These lists of graphs have subsequently been used to test conjectures and to search for potential 
counterexamples to others, as well as to make new discoveries, look for patterns, 
and formulate new conjectures, as for 
example in  \cite{DvDFGG, FengKwak, FengZhou, FerHuj, FreKut, FKPV, MonWeiss, Parker, PSV2}. 

\smallskip

\begin{figure}[hhh]
\includegraphics[width=0.85\textwidth]{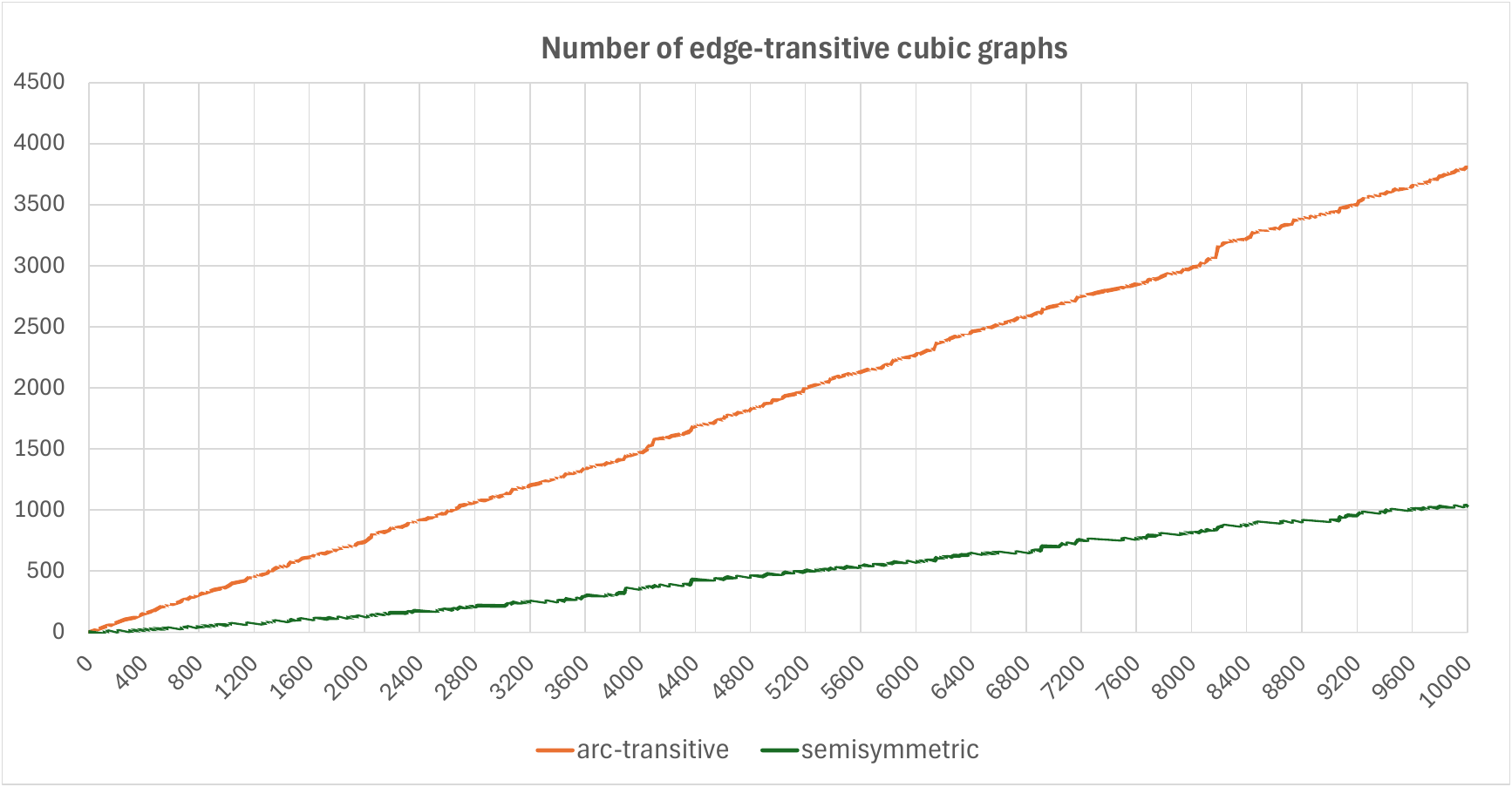} 
\caption{Numbers of edge-transitive cubic graphs 
of up to given order ($\le 10000$)} 
\label{fig:numET}
\end{figure}

One of the main aims of the current paper is to announce, document and explain the construction of a 
much larger census, namely one of all cubic edge-transitive on at most 10000 vertices. 
We can now state the following new theorem (accompanied by Figure~\ref{fig:numET} which illustrates 
the growth in the number of them),  while postponing a more detailed description 
until Section~\ref{sec:10000list}. 
   
 \begin{theorem}
 \label{the:census}
 There are $4858$ connected finite edge-transitive cubic graphs on up to $10000$ vertices, 
 with $3815$ of them being arc-transitive, and the other $1043$ being semisymmetric.
 \end{theorem}  
  
 \begin{remark}
 The $4858$ cubic edge-transitive graphs themselves are stored online
at \cite{OnlineCensus}, 
with summary data also available at \cite{ATConder,SemiSymConder}. 
At \cite{OnlineCensus}, for each admissible $n$ we provide files `CAT\_$n$.s6' and `CSS\_$n$.s6' containing the sparse6 codes \cite{sparse6}
 of cubic arc-transitive and cubic semisymmetric graphs of order $n$, respectively, with one graph per line.
For future reference, we denote the graph in the $k^{\rm th}$ row of the file CAT\_$n$.s6 as ${\rm CAT}(n,k)$ and
the one in the $k^{\rm th}$ row of the file CSS\_$n$.s6 as ${\rm CSS}(n,k)$.  Additionally, several  graph invariants have been pre-computed, and
are presented at \cite{OnlineCensus} in a tabular format. Moreover, hamiltonicity of each of these graphs  was tested using
the {\tt hamheuristic} command from the {\tt gtools} part of {\tt nauty} \cite{nauty}. Apart from two well-known exceptions (the Petersen graph CAT(10,1)
and the Coxeter graph CAT(28,1)),  each of the cubic edge-transitive graphs in the list was found to contain a hamilton cycle.  Hence in particular, all of them have a hamilton path.
\end{remark}
   
Let us now turn our attention to the only infinite graph that we will be interested in, 
namely the $3$-regular tree $\Tree_3$.
  
The automorphism group of $\Tree_3$, when equipped with the topology of point-wise convergence, is
 a locally compact totally disconnected topological group. 
 A subgroup $G$ of $\Aut(\Tree_3)$ is discrete in this topology if and only if the stabiliser in $G$ 
 of every vertex $v\in \V(\Tree_3)$ is finite.
By theorems of Tutte~\cite{T47}, Djokovi\'c and Miller~\cite{DjM} and Goldschmidt~\cite{Gold80}, 
there are precisely seven conjugacy classes of discrete arc-transitive subgroups of $\Aut(\Tree_3)$ 
and precisely fifteen conjugacy classes of discrete edge- but not vertex-transitive (and therefore semisymmetric) 
subgroups of $\Aut(\Tree_3)$. 

The seven classes of arc-transitive groups were determined and explicitly described in \cite{DjM}, 
with alternative descriptions as finitely-presented groups given in \cite{CondLor}, and will be denoted 
here by $\DjM_1$, $\DjM_2^{\,1}$, $\DjM_2^{\,2}$, $\DjM_3$, $\DjM_4^{\,1}$, $\DjM_4^{\,2}$ and $\DjM_5$.  
The fifteen classes of semisymmetric groups were determined in \cite{Gold80} and explicitly described 
as finitely-presented groups in \cite{CMMP}, and will be denoted here by $\Gold_1$, $\Gold_1^{\,1}$, $\Gold_1^{\,2}$, $\Gold_1^{\,3}$, $\Gold_2$, $\Gold_2^{\,1}$, $\Gold_2^{\,2}$, $\Gold_2^{\,3}$, $\Gold_2^{\,4}$, 
$\Gold_3$,    $\Gold_3^{\,1}$,    $\Gold_4$,    $\Gold_4^{\,1}$,    $\Gold_5$ and  $\Gold_5^{\,1}$.
(Incidentally, in the presentation for  $\Gold_3$ in \cite{CMMP}, the final relator should be $cdydy$, not $cdyd$. 
Curiously, it could also be replaced by $cd^{-1}ydy$, under an isomorphism that inverts each of the 
generators $c, d, x$ and $y$, which we will mention later in the paper.) 

\smallskip

Each of the $7+15 = 22$ classes of edge-transitive subgroups of $\Aut(\Tree_3)$ comes from an `amalgam'  
of finite groups, as we explain briefly below.

It follows from Bass-Serre theory that every discrete edge-transitive subgroup of $\Aut(\Tree_3)$ can be 
obtained as an amalgamated free product  $A*_C B$ of a triple of groups $(A,B,C)$ with $C=A\cap B$.  
The latter is called an {\em amalgam\/} of groups.  
In our context, for some connected finite cubic graph $\Gamma$, 
the following hold:
   \begin{enumerate}
   \item in the arc-transitive case:
   $A=G_v$, $C=G_{uv}$ and $B=G_{\{u,v\}}$ are the stabilisers of an \\ incident vertex-arc-edge triple in an arc-transitive group  $G$
   of automorphisms of $\Gamma$, while \\[-10pt]
   \item in the semisymmetric case:
   $A=G_v$, $C=G_{uv}$ and $B=G_u$ are the stabilisers of an \\ incident vertex-arc-vertex triple in a semisymmetric group $G$ of automorphisms of $\Gamma$.
   \end{enumerate}
  
In all of the $22$ amalgams associated with edge-transitive subgroups of $\Aut(\Tree_3)$, 
the groups $A$, $B$ and $C$ are finite.   
Moreover, no non-trivial subgroup of $C$ is normal in both $A$ and $B$, by faithfulness 
of the action of the group $G$ on the edges of the graph $\Gamma$, and also $(|A\!:\!C|, |B\!:\!C|)  = (3,3)$ or $(3,2)$, 
depending on whether the action of $G$ on $\Gamma$ is semisymmetric or arc-transitive. 
Amalgams satisfying these properties are said to be {\em finite\/} and {\em simple},
 with {\em index\/} $(3,3)$ or $(3,2)$ respectively.
   
Conversely, if $(A,B,C)$ is a finite simple amalgam with index $(3,3)$ or $(3,2)$, then the corresponding 
amalgamated product $A*_C B$ is obtainable from a representative of exactly one of the $22$ conjugacy 
classes of discrete edge-transitive subgroups of $\Aut(\Tree_3)$. 

In fact, this establishes a bijective correspondence  between the isomorphism classes of simple amalgams 
of index $(3,3)$ and $(3,2)$ and the $22$ conjugacy classes of edge-transitive discrete subgroups 
of $\Aut(\Tree_3)$, with the $15$ amalgams of index $(3,3)$ corresponding to semisymmetric subgroups, 
and the $7$ amalgams of index $(3,2)$ to arc-transitive subgroups.
(Here an {\em isomorphism\/} of two amalgams is to be understood in an appropriate category 
of diagrams $A \leftarrow C \rightarrow B$.)

\smallskip

We will sometimes abuse notation and use the above symbols for the $22$ classes of edge-transitive 
discrete subgroups of $\Aut(\Tree_3)$ to denote specific representatives of these classes, 
rather than the classes themselves, and for reasons that will become apparent later, 
we will also call the seven classes of discrete arc-transitive groups the {\em Djokovi\'c-Miller classes} (or {\em Djokovi\'c-Miller amalgams}), and the fifteen classes of discrete semisymmetric groups 
the {\em Goldschmidt classes} (or {\em Goldschmidt amalgams}). 
   
The importance for finite cubic edge-transitive graphs of the above classification of discrete edge-transitive 
subgroups of $\Aut(\Tree_3)$ is revealed by the following folklore theorem, which again follows from 
Bass-Serre theory. We state it in the context of cubic edge-transitive graphs, but of course it also holds 
in a more general setting.  

Also we note in advance that if $N$ is a normal subgroup of the class 
representative group $\tG$, then $\Tree_3/N$ is the {\em quotient graph}, whose vertices are the 
orbits of $N$ on $\V(\Tree_3)$, with two such orbits adjacent whenever they contain a pair of 
vertices that are adjacent in $\Tree_3$.

\begin{theorem}   
 \label{thm:serre}
Let $\Gamma$ be a connected finite cubic graph, and let $G$ be a group that acts transitively on the edges of $\Gamma$. 
Then there exists  a unique conjugacy class of discrete edge-transitive subgroups of $\Aut(\Tree_3)$ such that 
for one (and hence every) representative $\tG$ of that class, there exists a normal subgroup $N\norml \tG$  
with the following properties:
   \begin{enumerate}
    \item $N$ is a free group of finite rank $\beta$, where $\beta = |\V(\Gamma)| - |\E(\Gamma)| + 1$ is the Betti number of $\Gamma;$
    \item $N$ acts semiregularly on the vertices and on the edges of $\Tree_3$, 
    and the quotient graph $\Tree_3/N$ is isomorphic to $\Gamma;$
    \item the quotient group $\tG/N$ is isomorphic to $G$, and what is more, 
    any isomorphism between $\Gamma$ and $\Tree_3/N$, when taken together with an appropriate 
    isomorphism from $G$ to $\tG/N$, yields an isomorphism between the action of $G$ on $\V(\Gamma)$
    and the obvious action of $\tG/N$ on $\V(\Tree_3/N)$.
   \end{enumerate}
 \end{theorem}   
   
The above theorem reduces the problem of finding all connected cubic edge-transitive graphs 
with at most $m$ edges to the problem of finding all normal subgroups of index up to $c m$ 
in a representative $\tG$ from each of the Djokovi\'c-Miller classes and each of the Goldschmidt classes,
where $c$ is the order of the associated edge-stabiliser in $\tG$.

Of course the automorphism group of the graph $\Gamma \cong \Tree_3/N$ might be larger 
than $G \cong \tG/N$.  (In particular, if $G$ acts semisymmetrically on $\Gamma$, it can often 
happen that $\Gamma$ is vertex-transitive and hence arc-transitive.) 
In that case $\Aut(\Gamma)$ is a quotient of some amalgam `larger' than $\tG$, 
containing $\tG$ as a subgroup of finite index. 

Accordingly, although the $22$ classes of discrete edge-transitive of subgroups of $\Aut(\Tree_3)$ 
are well understood and descriptions of their representatives as finitely-presented groups are known, 
it is important to understand the possible inclusions of the associated amalgams in each other. 
Such inclusions are already known for the seven Djokovi\'c-Miller classes in the arc-transitive case, 
as explained in~\cite{DjM} and taken further in~\cite{CondLor} and~\cite{CN}, and we complement 
that work by dealing with inclusions of each of  the $15$ Goldschmidt classes in one or more of the other $21$. 

\begin{figure}[hhh]
\phantom{}
\includegraphics[width=1.02\textwidth]{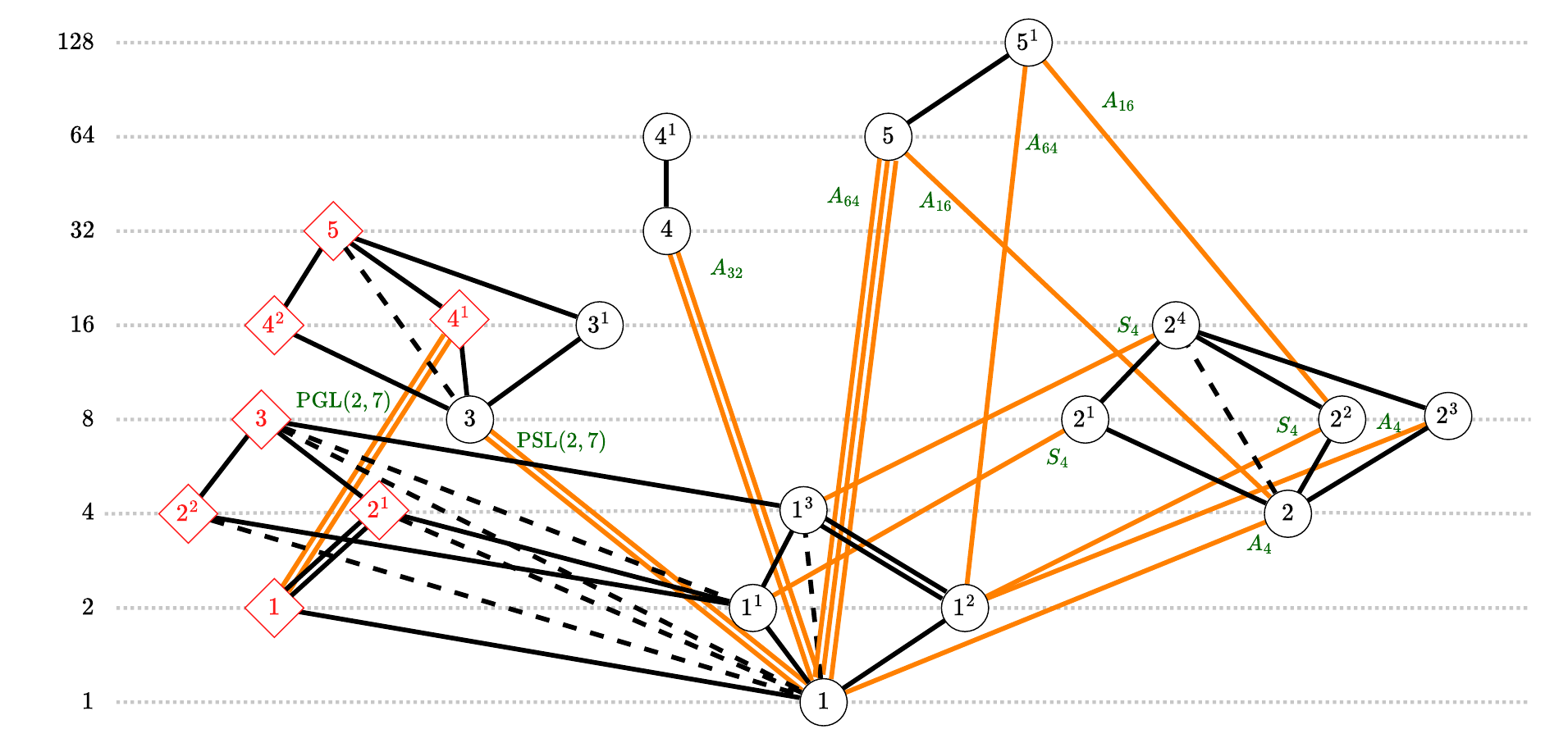} 
\caption{
The $22$ classes of discrete edge-transitive subgroups of $\Aut(\Tree_3)$.} 
\label{fig:GoldDjM}
\end{figure}

To be more precise, suppose that $\cC_1$ and $\cC_2$ are two of the $22$ classes. 
If a representative of $\cC_1$ is a subgroup of a representative of $\cC_2$, then we say 
that $\cC_1$ is {\em included} in $\cC_2$, and write $\cC_1 \le \cC_2$. 
Moreover, if $\cP$ is a property that is meaningful for a pair of a groups, one contained in the other, such as `normal', or `non-normal', or `maximal', or `non-maximal',
 then we say that $\cC_1$ is $\cP$-included  in $\cC_2$ provided that representatives $X\in \cC_1$ and $Y\in \cC_2$ can be chosen so that
the pair $X\le Y$ has property $\cP$. 
It is not entirely obvious, but we will indirectly show that if $\cC_1$ is normally included in $\cC_2$,
then it cannot be non-normally included in $\cC_2$.    There are cases, however, 
in which $\cC_1$ is both maximally and non-maximally included in $\cC_2$.
For example, a representative $G$ of $\Gold_5$ contains representatives of $\Gold_1$ as maximal subgroups, but also other representatives of $\Gold_1$ that are not maximal in $G$, being contained in a representative  
of $\Gold_2$, which in turn is maximal in $G$. 
Similarly, representatives of $\Gold_1^{\,2}$ are included in representatives of $\Gold_5^{\,1}$ 
both maximally as well as
non-maximally via $\Gold_2^{\,2}$.  
Moreover, as can be seen from Figure~{\em\ref{fig:GoldDjM}}, the are no other examples of 
simultaneous maximal and non-maximal inclusions among Goldschmidt and Djokovi\'{c}-Miller classes. 
Details can be found in Subsection~\ref{sec:InclusionsGtoG}.

We can now state the second of our main new theorems:
   
\begin{theorem}
\label{thm:Inclusions}
 The mutual inclusions of the $22$ conjugacy classes of discrete edge-transitive subgroup of $\Aut(\Tree_3)$ 
 are as indicated in Figure~{\em\ref{fig:GoldDjM}}.
 \end{theorem}

In Figure~\ref{fig:GoldDjM}, the classes are arranged in horizontal layers, with each layer corresponding 
to the order (indicated by the number left of each layer) of the edge-stabiliser in a representative of the class.

The seven Djokovi\'c-Miller classes are represented by diamond shapes of red colour, while the fifteen Goldschmidt classes are encircled and written in black.
Solid lines represent `maximal inclusions' (inclusions of a representative subgroup as a maximal subgroup 
of a representative subgroup of another). In particular, black solid lines represent normal maximal inclusions, 
while orange solid lines represent non-normal maximal inclusions.  
Of course all index $2$ subgroups are maximal.

We have also indicated all normal inclusions that are not maximal: they are depicted by dashed black lines.
The double black line between the Djokovi\'c-Miller classes $\DjM_1$ and $\DjM_2^{\,1}$ and between Goldschmidt  classes
$\Gold_1^{\,2}$ and $\Gold_1^{\,3}$ means that the larger group of one type  contains two smaller normal subgroups of the other type.
Similarly, double and triple orange lines (such as between $\DjM_1$ and $\DjM_4^{\,1}$, or between $\Gold_1$ and $\Gold_5$) indicate that the larger group contains two or three conjugacy classes of maximal
subgroups of a given type.

The green text next to the orange lines gives information about the quotient of the top group, say $K$, by the core of the (non-normal) subgroup, say $J$. 
When the core is the alternating group $A_n$ or symmetric group $S_n$, then the conjugacy class always has size $n$.
On the other hand, the group in the case of the group $\PSL(2,7)$, the size of the conjugacy class is $8$ and the action on the corresponding coset space
 is as on the projective line $\PG(1,7)$.

More details of the inclusions and further explanations are given in Section~\ref{sec:Inclusions}.

\medskip

A quick look at Figure~\ref{fig:GoldDjM} reveals  that among the Goldschmidt classes, only the classes $\Gold_1$, $\Gold_1^{\,j}$ for $j\in  \{1,\ldots, 3\}$,
$\Gold_3$ and $\Gold_3^{\,1}$ are included in any of the Djokovi\'c-Miller classes. This has the following immediate consequence:

\begin{corollary}
 If a finite cubic graph admits a semisymmetric group of automorphism of one of the following types: $\Gold_2$, $\Gold_2^{\,j}$ for some $j\in \{1,\ldots,4\}$,
 $\Gold_4$,  $\Gold_4^{\,1}$,  $\Gold_5$ or  $\Gold_5^{\,1}$, then the graph is semisymmetric.
\end{corollary}

 Let us now turn our attention back to finite graphs. As we  pointed out in Theorem~\ref{thm:serre}, each group $G$ acting edge-transitively on
 a finite cubic graph is a quotient of a representative of exactly one of the $22$ Djokovi\'c-Miller and Goldschmidt classes, and we will then say
 that this class is the {\em type} of $G$. In fact, each class $\cC$ can be realised in a finite edge-transitive cubic graph, in the sense that there exists such a graph $\Gamma$ and an edge-transitive group $G$ of automorphisms of $\Gamma$ such that $G$ is of type $\cC$.
 
 This can be seen directly, by exhibiting a carefully chosen sample of finite cubic graphs whose automorphism groups contain edge-transitive subgroups.
 In particular, the Djokovi\'c-Miller classes $\DjM_1$, $\DjM_2^{\,1}$, $\DjM_2^{\,2}$ and $\DjM_2^{\,3}$ occur in the complete bipartite graph $\rK_{3,3}$,
  denoted by $\CAT(6,1)$ in our census, while the types $\DjM_4^{\,1}$, $\DjM_4^{\,2}$ and $\DjM_5$ occur in Tutte's 8-cage on $30$ vertices, denoted by $\CAT(30,1)$.
Moreover, $\rK_{3,3}$ and Tutte's $8$-cage contain semisymmetric groups of types $\Gold_1$, $\Gold_1^{\,1}$, $\Gold_1^{\,2}$, $\Gold_1^{\,3}$, and
$\Gold_3^{\,1}$, $\Gold_3^{\,2}$, respectively. Similarly, types $\Gold_2$, $\Gold_2^{\,1}$, $\Gold_2^{\,2}$, $\Gold_2^{\,3}$ and $\Gold_2^{\,4}$ are realised in the
smallest cubic semisymmetric graph, namely Marion Gray's graph on $54$ vertices, denoted by $\CSS(54,1)$, while the types $\Gold_4$ and $\Gold_4^{\,1}$
occur in the incidence graph of the generalised hexagon, denoted by $\CSS(126,1)$, and the types $\Gold_5$ and $\Gold_5^{\,1}$ occur in
the unique semisymmetric graph on $990$ vertices (of girth $16$ and diameter $12$), denoted by $\CSS(990,1)$.
 
Alternatively, one can prove the same fact using a theorem of Baumslag~\cite[Theorem 2]{Baum}, which states that an amalgamated free product of
 finite groups is residually finite.
 
It is not so immediately clear, however, that (or how) each of the $22$ classes can be realised as the {\em full\/}  automorphism group of a cubic edge-transitive graph.

Let us say that a conjugacy class $\cC$ of discrete edge-transitive subgroup of $\Aut(\Tree_3)$ has a {\em strong realisation} in a graph $\Gamma$, provided that $\Aut(\Gamma)$ is of type $\cC$. When that happens, we shall say that the graph $\Gamma$ has {\em type} $\cC$.
As we will show in this paper, all Djokovi\'c-Miller and Goldschmidt classes can indeed be realised strongly.
 We can prove this is in two different ways:  first constructively, by
 finding examples among the graphs of order at most 10000, or by ad-hoc constructions when no such small examples exist,
 and then also non-constructively, using the theory of lifting groups along regular covering projections of graphs.
  The latter approach proves for each type $\cC$ the existence of not just one but an infinite family of examples of graphs strongly realising $\cC$.
   To summarise, we prove the following:
  
\begin{theorem} 
\label{thm:existence}
For each conjugacy class $\cC$ of discrete edge-transitive subgroups of $\Aut(\Tree_3)$, 
there exists an infinite family of finite cubic edge-transitive graphs whose full automorphism 
groups are of type~$\cC$.
\end{theorem}
  
Existence of an infinite family of graphs strongly realising each Djokovi\'c-Miller or Goldschmidt class 
raises a question about the growth of the number of such graphs with respect to their order. 
For a conjugacy class $\cC$ of discrete edge-transitive subgroup of $\Aut(\Tree_3)$, 
let $f_\cC(n)$ denote the number of cubic edge-transitive graphs of type $\cC$ on at most $n$ vertices. 
The sums of all $f_\cC(n)$ over all Djokovi\'c-Miller class $\cC$, as well as the sum over all Goldschmidt classes,  are shown in Figure~\ref{fig:numET} for $n\le 10000$. 

\begin{figure}[hhh]
\begin{tabular}{cc}
\hspace{-3mm}
\includegraphics[width=0.478\textwidth]{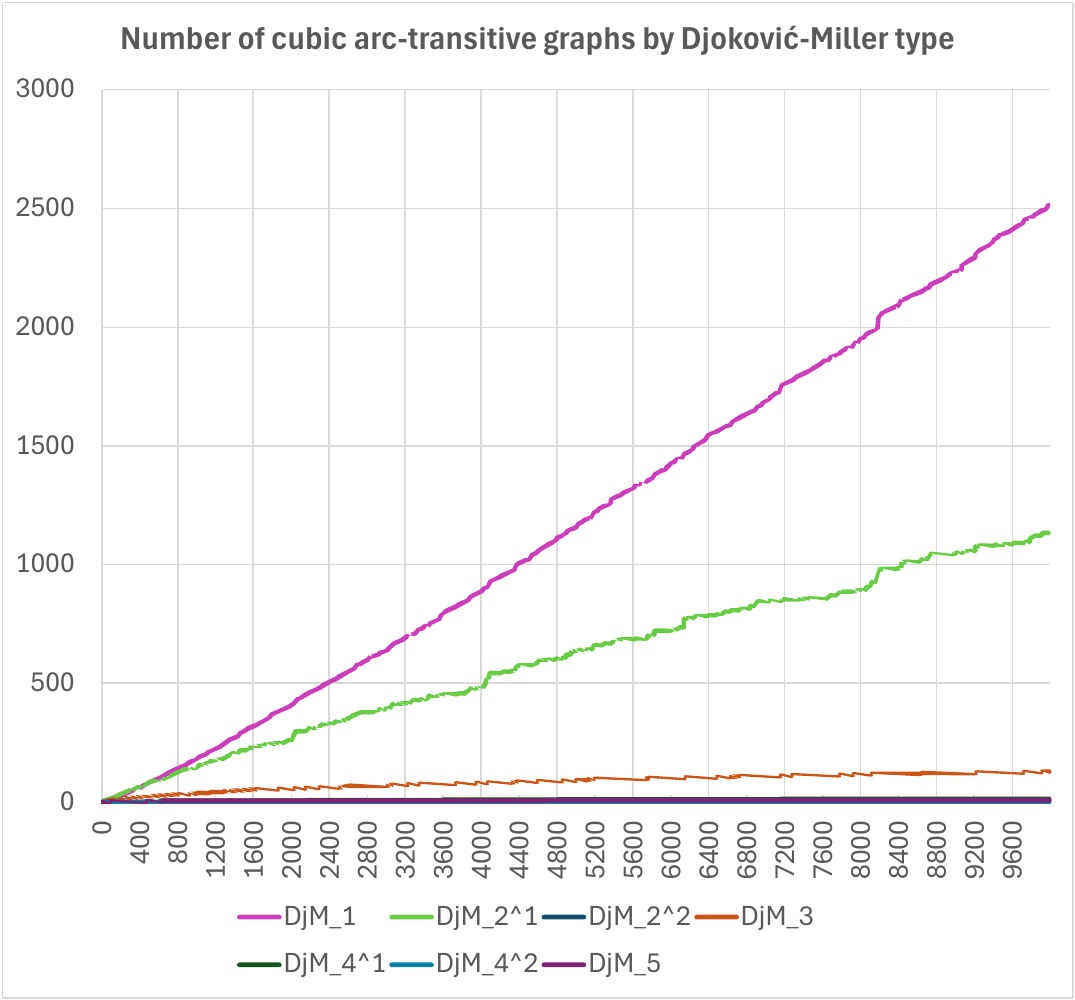} &
\includegraphics[width=0.5\textwidth]{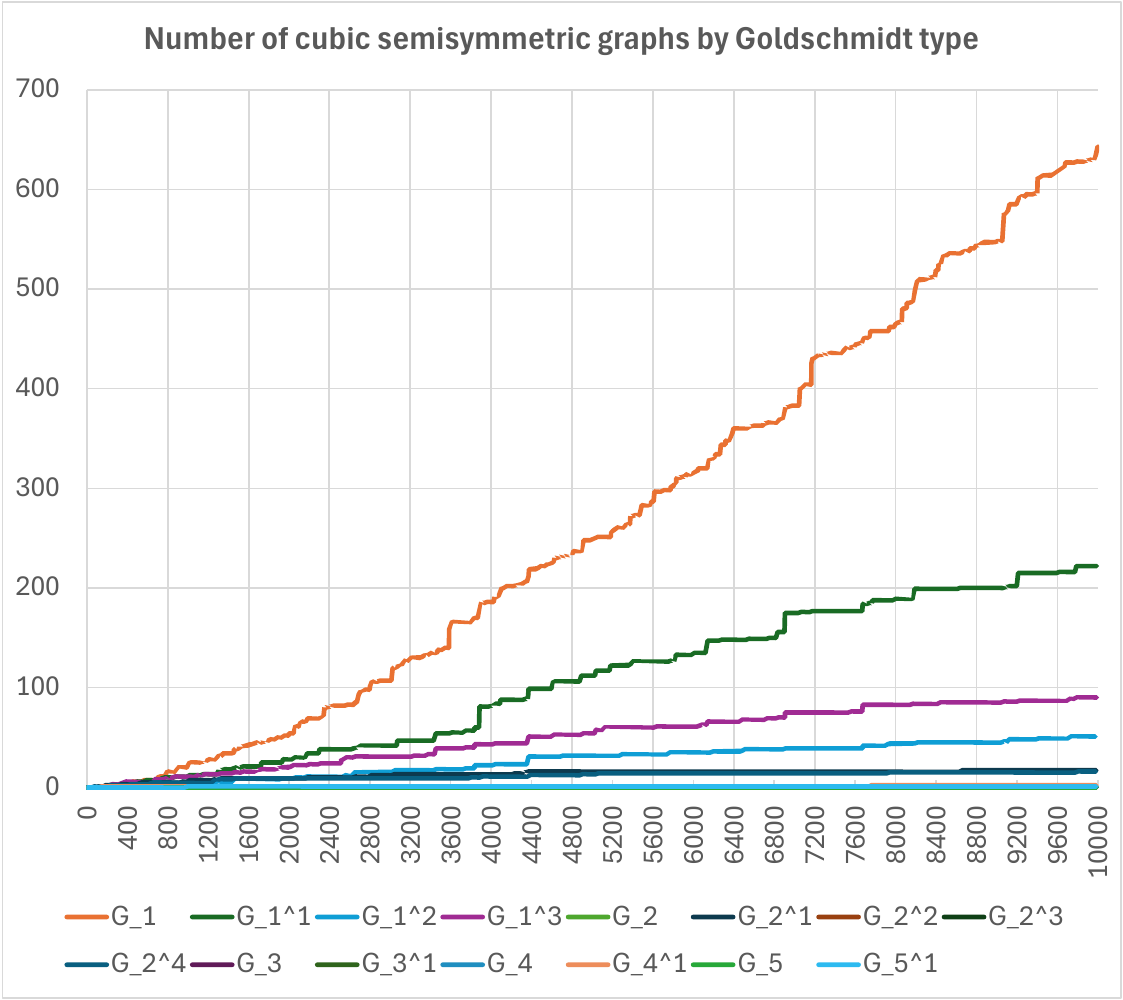}
\end{tabular}
\caption{%
Numbers of arc-transitive and semisymmetric cubic graphs of each type. 
}
\label{fig:numClasses}
\end{figure}

A quick look might suggest that the growth of these sums is roughly linear. 
 A similar conclusion might be drawn from the plots of each of the corresponding $22$ functions $f_\cC(n)$; see Figure~\ref{fig:numClasses}.

The lines on the left-hand part of Figure~\ref{fig:numClasses}, from the top to bottom, 
represent the classes $\DjM_1$, $\DjM_2^{\,1}$, $\DjM_3$, $\DjM_4^{\,1}$, $\DjM_2^{\,2}$, $\DjM_5$ and $\DjM_4^{\,2}$, with the last four essentially indistinguishable. On the right-hand part, the Goldschmidt classes 
appear in the order $\Gold_1$, $\Gold_1^{\,1}$, $\Gold_1^{\,3}$, 
$\Gold_1^{\,2}$, then $\Gold_2^{\,1}$ and $\Gold_2^{\,4}$ almost overlapping, similarly, $\Gold_4^{\,1}$, $\Gold_3^{\,1}$, and $\Gold_5^{\,1}$ lying just above the horizontal axis,
and the lines corresponding to $\Gold_2$, $\Gold_2^{\,2}$, $\Gold_2^{\,3}$, $\Gold_3$, $\Gold_4$ and $\Gold_5$ coinciding with the horizontal axis (since there are no graphs of these types having order at most $10000$).

But such conclusions about linear growth of the sums of all $f_\cC(n)$ for a given class $\cC$ 
would be far from true, because we prove the following in Section~\ref{sec:Enumeration}:

\begin{theorem}
\label{thm:enum}
For each conjugacy class $\cC$ of discrete edge-transitive subgroups of $\Aut(\Tree_3)$, there exist
positive real constants $a$ and $b$ and a positive integer $n_0$ such that 
$$
   n^{a \log(n)} \le f_\cC(n) \le    n^{b \log(n)} \quad \hbox{for all} \ n\ge n_0.
$$
\end{theorem}

This asymptotic enumeration theorem counters not only an impression that one might gain by considering graphs of orders up to $10000$,
but also a widely spread belief that cubic semisymmetric graphs (or even those of a given class) are rare.

\begin{figure}[hhh]
\includegraphics[width=0.65 \textwidth]{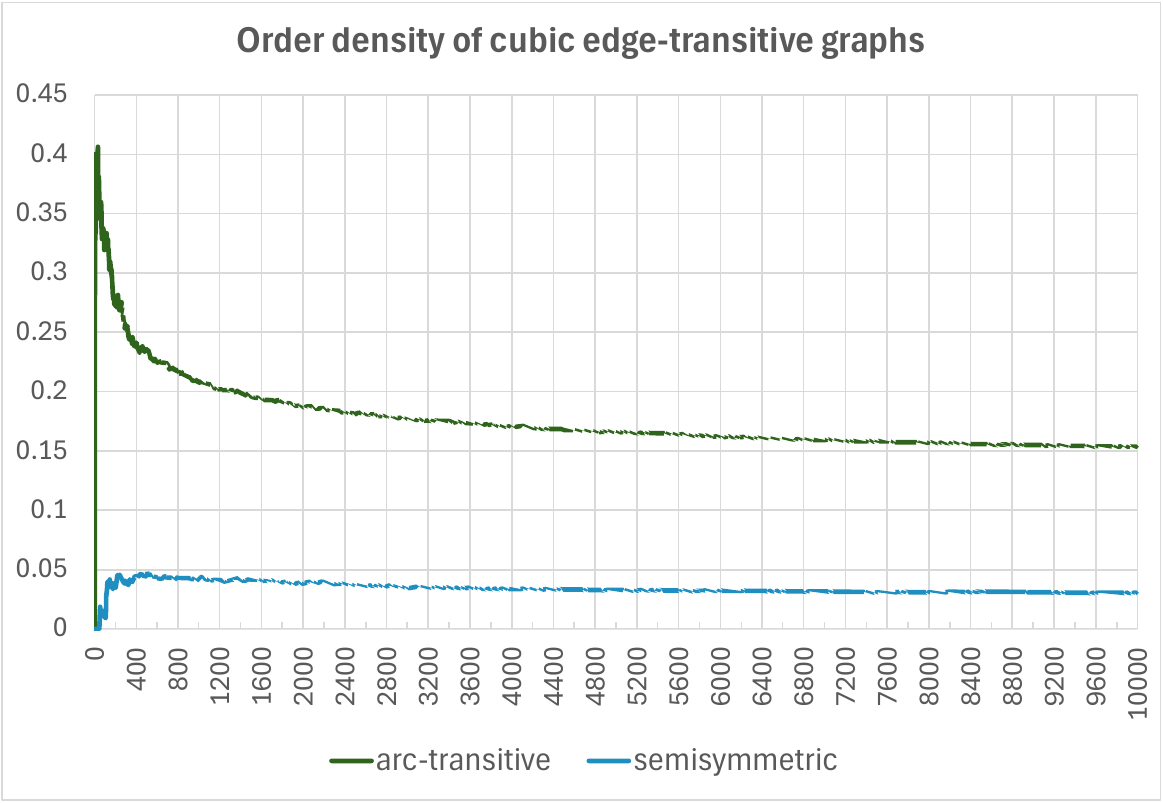}
\caption{%
The order density of cubic arc-transitive and semisymmetric graphs.
}
\label{fig:density}
\end{figure}

The fact that the number of cubic edge-transitive graphs grows faster than any polynomial function might also suggest that the set of their orders represents a significant proportion of the positive integers. 
This expectation is to some extent supported by the
empirical data that can be derived from our lists of cubic edge-transitive graphs, as depicted in Figure~\ref{fig:density}. 
This figure exhibits the proportion of positive integers $n$ that are the orders of arc-transitive cubic graphs, 
and the same for semisymmetric cubic graphs, which we call the {\em order density} in each case.

Looking at Figure~\ref{fig:density}, one might conjecture that this ratio stabilises at about $15\%$ for arc-transitive graphs, and at about $3\%$ for semisymmetric graphs. 
Such a conjecture would also be wrong, however, because the natural density of the set of orders of connected finite cubic edge-transitive graphs is $0$,
as was proved recently by Conder, Verret and Young (see \cite[Theorem 4.3]{CVY}).

\smallskip

We now return to the matter of determining edge-transitive cubic graphs of up to a given order. 
 
In the arc-transitive case, we have already mentioned the Foster census \cite{Foster} (for order up to $512$) 
and its completion and extension up to order $768$ in \cite{CD} and later to order $2048$.   
The two extensions were obtained by finding all suitable quotients of each of the seven Djokovi\'c-Miller groups. 
The construction of the list of semisymmetric cubic graphs on up to $768$ vertices in \cite{CMMP} 
was achieved using Goldschmidt's classification \cite{Gold80} of finite primitive amalgams of index $(3,3)$. 

The Goldschmidt classification implies that 
every edge- but not vertex-transitive group of automorphisms of a connected finite 
cubic graph is obtainable as a homomorphic image of the universal completion of one of the $15$ Goldschmidt  
amalgams we described earlier.  
In all $15$ of these finite primitive amalgams of index $(3,3)$, the orders of the two subgroups 
$A$ ($\cong G_u$) and $B$ ($\cong G_v$) are at most $384 = 3 \cdot 2^7$, 
and it follows that in the automorphism group of any connected finite semisymmetric cubic graph $\Gamma$, 
the stabiliser of any given vertex has order $3\cdot 2^{s-1}$ for some $s\le 8$. 
Clearly this is the analogue for semisymmetric cubic graphs of the now classical theorem of Tutte \cite{T47}, 
which states that $s\le 5$ in the case where the graph $\Gamma$ is arc-transitive, 
and in particular, the word `semisymmetric' can be replaced by `edge-transitive' in the previous sentence. 

\smallskip

This paper is organised as follows.  

In Section~\ref{sec:Inclusions} we consider the fifteen Goldschmidt 
classes and the seven Djokovi\'c-Miller classes in more detail, determining all possible inclusions 
of such amalgams among each other.  Then in Section~\ref{sec:10000list} we describe how we used 
the Djokovi\'c-Miller and Goldschmidt amalgams to help find all the edge-transitive cubic graphs 
of order up to 10000, and in Section~\ref{sec:Examples} we give some information about examples 
of connected finite edge-transitive cubic graphs of each type, including the smallest one in all but a few cases. 
Then we use regular covering projections in Section~\ref{sec:Existence} to prove our Theorem~\ref{thm:existence} about there being infinitely many examples of each type, 
and finally we justify our comments above about asymptotic enumeration of edge-transitive graphs 
of each type in Section~\ref{sec:Enumeration}. 

Just before continuing, however, we comment on notation used in various parts of the paper. 
We use $C_n$ to denote the cyclic group of order $n$ (or degree $n$ when considered as a permutation group), 
and $D_n$, $A_n$ and $S_n$ to denote the dihedral, alternating and symmetric groups of degree $n$. 
(This contrasts with the notation in \cite{CMMP}, where $D_{2n}$ was used for the dihedral group of order $2n$.)

\section{Conjugacy classes of discrete edge-transitive subgroups of $\Aut(\Tree_3)$}
\label{sec:Inclusions}

Recall that the fifteen isomorphism classes of finite simple amalgams of index $(3,3)$ 
were determined by Goldschmidt in \cite{Gold80}, and explicit descriptions of the corresponding 
semisymmetric subgroups of $\Aut(\Tree_3)$ as finitely-presented groups were given in \cite{CMMP}.  
Similarly, the seven isomorphism classes of finite simple amalgams of index $(3,2)$ 
were determined by Djokovi\'c and Miller in \cite{Gold80} in \cite{DjM}. 

In the rest of this section, we will consider each of the Goldschmidt and Djokovi\'c-Miller classes 
of subgroups of $\Aut(\Tree_3)$, one by one, determining inclusions of other (`smaller') amalgams 
within them. This was already done in \cite{DjM}, \cite{CondLor} and \cite{CN} for the inclusions among 
the seven Djokovi\'c-Miller classes, but we repeat the information here for the sake of completeness. 
We consider the inclusions of Goldschmidt amalgams among each other in Subsection \ref{sec:InclusionsGtoG}, 
and the inclusion of both kinds of amalgams in Djokovi\'c-Miller amalgams 
in Subsection~\ref{sec:InclusionsGtoDM}.

\subsection{Inclusions among Goldschmidt classes}
\label{sec:InclusionsGtoG}

${}$ 

\par 

To find out which of the Goldschmidt classes are included in other such classes, 
we made extensive use of the {\sc Magma} system~\cite{Magma}. 

To explain this in more detail, we suppose that we fix a Goldschmidt class ${\cal G}$ and a semisymmetric subgroup $G\le \Aut(\Tree_3)$ representing it.
Then $G=A*_C B$, where $A=G_v$, $C=G_{uv}$ and $B=G_u$ for two adjacent vertices $v$ and $u$ of $\Tree_3$, and then look inside $G$ for an edge-transitive subgroup $H$ from some  `smaller' Goldschmidt class ${\cal H}$. 
Note that the amalgam $(H_v,H_{uv},H_u)$ is then a representative of the class ${\cal H}$, 
and that edge-transitivity of $H$  implies that $G=G_{uv} H$, and hence that $|G\!:\!H| =|G_{uv}| / |H_{uv}|$.

In most cases (and whenever $|G_{uv}| \le 16$), we used the {\tt LowIndexSubgroups} procedure to find all conjugacy classes 
of subgroups of the relevant index (such as index $|G_{uv}|/2 = 4$ when determining edge-transitive subgroups of a group $G$ of type $\Gold_2^{\,1}$ 
that belong to the class $\Gold_1^{\,1}$). Once all conjugacy classes of subgroups of up to desired index are found, edge-transitivity of such a subgroup $H$
could be tested by checking whether $G_{uv}$ acts transitively on the coset space $(G\!:\!H)$.

When the relevant index was $32$ or more, the latter approach was too slow, and so instead we needed 
to take a different approach, namely as follows.  

If $(u,v)$ is a given arc, and the automorphisms $x$ and $y$ of order $3$ in $\Aut(\Tree_3)$ 
are chosen to fix the vertices $u$ and $v$ respectively, then any copy of $\Gold_1$ is conjugate to 
a subgroup generated by two elements $x'$ and $y'$  of order 3 fixing $u$ and $v$ (respectively), 
such that $x'x^{-1}$ and $y'y^{-1}$ fix the arc $(u,v)$. 
Here both $x'x^{-1}$ and $y'y^{-1}$ lie in the stabiliser $C$ of $(u,v)$, making $x' = rx$ and $y' = sy$ for some $r$ and $s$  in $C$, and as $C$ is small, it is easy to run through the few possibilities for the pair $(r,s)$, and for each one,  test the subgroup generated by $x' = rx$ and $y' = sy$. 
This helps to find all possibilities for a conjugacy class of edge-transitive subgroups of type $\Gold_1$ 
with the required index in $G$, and can be adapted to find edge-transitive subgroups of other types 
(with smaller index in $G$) as necessary.

\smallskip

In what follows, we consider each of the fifteen Goldschmidt classes in turn. 
We first describe each class $\cG$ by specifying a finite presentation for a representative $G \in \cG$, 
and for the purposes of this exposition, we will call this presentation a {\em standard presentation} for the 
corresponding Goldschmidt class, and call the generators in this presentation the {\em standard generators} 
for that class. 
We also list generators for the stabilisers $A=G_u$, $B=G_v$ and $C= A \cap B = G_{uv}$. 
We then list all edge-transitive subgroups of $G$ (one per each $G$-conjugacy class) and provide information about the sizes of their $G$-conjugacy classes and
their Goldschmidt types. Note that if a subgroup $H\le G$ is of Goldschmidt type $\cH$, then it can be presented in terms of the standard presentation of the class $\cH$. 
 The inclusion of $H$ to $G$ is then described by specifying the images of the standard generators of $\cH$ 
 as words in the standard generators of $\cG$. 
 
 For example, if the standard presentations of $\cH$ and $\cG$ are 
  $H=\langle\, c, x, y \mid c^2, x^3, y^3, (cx)^2, (cy)^2 \,\rangle$ 
  and $G=\langle\, c, d, x, y \mid c^2, d^{2}, [c,d], x^3, y^3, (cx)^2, [d,x], [c,y], (dy)^2 \,\rangle$, 
then by $ (c,x,y)_H \mapsto (cd,x,y)_G$ we mean a homomorphism from $H$ into $G$ that maps the generators $c,x$ and $y$ of $H$ to the elements
 represented by  the words $cd,x$ and $y$ in the generators $c, d, x$ and $y$ of $G$. 
 
 Finally, for each Goldschmidt class we also provide a parameter called its {\em local arc-transitivity}, 
 defined as follows.
 For a positive integer $s$, an {\em $s$-arc} in a graph $\Gamma$ is a sequence of $s$ vertices, where any two consecutive vertices are adjacent and any three are pairwise distinct.
 Clearly, if $G\le \Aut(\Gamma)$, then the vertex-stabiliser $G_u$ acts on the set of all $s$-arcs
 whose initial vertex is $u$, in an obvious way. If this action is transitive, then we say that $G$ is locally $s$-arc-transitive at $u$, and the maximum value of $s$ such that $G$ is locally $s$-arc-transitive at
 $u$ is denoted by $s_u(G)$. If $G$ is an edge-transitive group of automorphisms of $\Gamma$ and $(u,v)$ is an arc of $\Gamma$, then the ordered pair $(s_u(G),s_v(G))$ is called the {\em local $s$-arc-transitivity} of $G$.
(Note that  this parameter does not depend on the choice of the given edge $\{u,v\}$.) The values $s_u(G)$ and $s_v(G)$ are both given, with the former referring to the vertex $u$ (whose stabiliser $G_u$ is the group $A$), and the latter referring to the vertex $v$ (whose stabiliser $G_v$ is the group $B$).
Of course the roles of $u$ and $v$ can be interchanged (in a form of duality), provided that the roles of $A$ 
and $B$ (and their generating sets) are interchanged at the same time.
Also to avoid possible notational confusion, note that in the case of $\Gold_5^{\,1}$, the vertex $v$ is not the same 
as the generator $v$ used in the presentation for the representative group $G$.
\bigskip


\noindent
\begin{tabular}{|l|}
\hline\hline
\\[-10pt] 
$\TypeSemicolon \Gold_1,\>\>\>  G=\langle\, x, y \mid x^3, y^3 \,\rangle$     \\
\\[-10pt] 
\hline
\\[-10pt] 
 $A = \langle x \rangle \cong C_3$ and 
 $B = \langle y \rangle \cong C_3$, of order $3$; \ 
 $C = \{1\}$ (trivial), $\>\> (s_u,s_v) = (1,1)$ 
\\[-10pt]
${}$ \\
\hline\hline
\end{tabular}
\bigskip

A group of type $\Gold_1$ contains no proper edge-transitive subgroup of $\Aut(\Tree_3)$.

\bigskip


\noindent 
\begin{tabular}{|l|}
\hline
\hline
\\[-10pt] 
$\TypeSemicolon  \Gold_1^{\,1},\>\>\> G = \langle\, c, x, y \mid c^2, x^3, y^3, (cx)^2, (cy)^2 \,\rangle$     \\
\\[-10pt] 
\hline
\\[-10pt] 
 $A = \langle c, x \rangle \cong S_3$ and 
 $B = \langle c, y \rangle \cong  S_3$, of order $6$; \ 
 $C =  \langle c \rangle \cong C_2$, $\>\> (s_u,s_v) = (2,2)$ 
\\[-10pt] 
${}$ \\
\hline\hline
\end{tabular}
\bigskip

\noindent
List of edge-transitive subgroups of $G$:
\medskip

\begin{itemize}[leftmargin=0.60cm] 
\item[$\bullet$] Index $2$ normal subgroup $H$ of type   $\Gold_1$, with standard generator 
   inclusion $(x, y)_H \mapsto (x, y)_G$, \\ and conjugation by $c$ in  $G$ inducing
   an automorphism of  $H$ that takes $(x,y)_H$ to $(x^{-1},y^{-1})_H$.
\end{itemize}
\bigskip


\noindent 
\begin{tabular}{|l|}
\hline
\hline
\\[-10pt] 
$\TypeSemicolon  \Gold_1^{\,2},\>\>\>  G = \langle\, c, x, y \mid c^2, x^3, y^3, (cx)^2, [c,y] \,\rangle$     \\
\\[-10pt] 
\hline
\\[-10pt] 
 $A = \langle c, x \rangle \cong S_3$ and 
 $B = \langle c, y \rangle \cong  C_6$, of order $6$; \ 
 $C =  \langle c \rangle \cong C_2$, $\>\> (s_u,s_v) = (1,2)$
 \\[-10pt] 
${}$ \\
\hline\hline
\end{tabular}
\bigskip

\noindent
List of edge-transitive subgroups of $G$:
\medskip

\begin{itemize}[leftmargin=0.60cm] 
\item Index $2$ normal subgroup $H$ of type   $\Gold_1$, with standard generator
   inclusion $(x, y)_H \mapsto (x, y)_G$, \\ and conjugation by $c$ in  $G$ inducing 
   an automorphism of  $H$ that takes $(x,y)_H$ to $(x^{-1},y)_H$. 
\end{itemize}
\bigskip


\noindent 
\begin{tabular}{|l|}
\hline
\hline
\\[-10pt] 
$\TypeSemicolon  \Gold_1^{\,3},\>\>\>  G = \langle\, c, d, x, y \mid c^2, d^2, [c,d], x^3, y^3, (cx)^2, [d,x], [c,y], (dy)^2 \,\rangle$     \\
\\[-10pt] 
\hline
\\[-10pt] 
 $A = \langle c, d, x  \rangle \cong D_6$ and 
 $B = \langle c, d, y \rangle \cong D_6$, of order $12$;  \ 
$C =  \langle c, d \rangle \cong C_2 \times C_2$, $\>\> (s_u,s_v) = (3,3)$
 \\[-10pt] 
${}$ \\
\hline\hline
\end{tabular}
\bigskip

\noindent
List of edge-transitive subgroups of $G$:
\medskip

\begin{itemize}[leftmargin=0.60cm]\setlength{\itemsep}{8pt}
\item Index $4$ normal (and non-maximal) subgroup $H$ of type   $\Gold_1$, with standard generator
   inclusion $(x, y)_H \mapsto (x, y)_G$, and conjugation by $c$ and $d$ in  $G$ 
   inducing automorphisms of  $H$ that take \\ $(x,y)_H$ to $(x^{-1},y)_H$ and $(x,y^{-1})_H$. 
\item Index  $2$ normal subgroup $H$ of type   $\Gold_1^{\,1}$, with standard generator
   inclusion $(c, x, y)_H \mapsto (cd, x, y)_G$, \\ and conjugation by $d$ in  $G$
   inducing an automorphism of  $H$ that takes $(cd,x,y)_H$ to $(cd,x,y^{-1})_H$. 
\item Index $2$ normal subgroup $H$ of type   $\Gold_1^{\,2}$, with standard generator 
   inclusion $(c, x, y)_H \mapsto (c, x, y)_G$, \\ and conjugation by $d$ in  $G$
   inducing an automorphism of  $H$ that takes $(c,x,y)_H$ to $(c,x,y^{-1})_H$. 
\item Another index $2$ normal subgroup $H$ of type   $\Gold_1^{\,2}$, with standard generator 
   inclusion given by $(c, x, y)_H \mapsto (d, y, x)_G$, and conjugation by $c$ in  $G$
   inducing an automorphism of  $H$ that takes $(c,x,y)_H$ to $(c,x^{-1},y)_H$. 
\end{itemize}

\smallskip
\noindent
Note: The last two subgroups (of type $\Gold_1^{\,2}$) are conjugate to each other 
within a group of type $\DjM_3$.

\smallskip\bigskip


\noindent 
\begin{tabular}{|l|}
\hline
\hline
\\[-10pt] 
$\TypeSemicolon  \Gold_2,\>\>\>  G = \langle\, c, d, x, y \mid 
c^2, d^2, [c,d], x^3, y^3, (cx)^2, [d,x], (cy)^3, dy^{-1}cy \,\rangle$ \\
\\[-10pt] 
\hline
\\[-10pt] 
 $A =  \langle c,d,x \rangle \cong D_6$ and 
 $B = \langle c,d,y \rangle \cong A_4$, of order $12$; \ 
 $C = \langle c,d  \rangle \cong C_2 \times C_2$, $\>\> (s_u,s_v) = (1,2)$
 \\[-10pt] 
${}$ \\
\hline\hline
\end{tabular}
\bigskip

\noindent
List of edge-transitive subgroups of $G$:
\medskip

\begin{itemize}[leftmargin=0.60cm]\setlength{\itemsep}{8pt}
\item Conjugacy class of four index $4$ maximal subgroups of type   $\Gold_1$, with
   standard generator inclusion for one of them (say $H$) given  by $(x, y)_H \mapsto (x, y)_G$, and quotient
   of  $G$ by the core of each subgroup being isomorphic to $A_4$.
\end{itemize}
\bigskip


\noindent 
\begin{tabular}{|l|}
\hline
\hline
\\[-10pt] 
$\TypeSemicolon  \Gold_2^{\,1},\>\>\>  G = \langle\, c, d, x, y \mid 
c^2, d^4, (cd)^2, x^3, (cx)^2, [d,x], y^3, (dy^{-1})^2, cydy \,\rangle$ \\
\\[-10pt] 
\hline
\\[-10pt] 
 $A =  \langle c,d,x \rangle \cong D_{12}$ and 
 $B = \langle c,d,y \rangle \cong S_4$, of order $24$; \ 
 $C = \langle c,d  \rangle \cong D_4$, $\>\> (s_u,s_v) = (3,4)$
 \\[-10pt] 
${}$ \\
\hline\hline
\end{tabular}
\bigskip

\noindent
List of edge-transitive subgroups of $G$:
\medskip

\begin{itemize}[leftmargin=0.60cm]\setlength{\itemsep}{8pt}
\item Conjugacy class of four index $8$ non-normal (and non-maximal) subgroups of type   $\Gold_1$, 
   with standard generator inclusion for one of them (say $H$) given  by $(x, y)_H \mapsto (x, y)_G$, and 
   quotient of  $G$ by the core of each subgroup being isomorphic to $S_4$.
\item Conjugacy class of four index $4$ non-normal  maximal subgroups of type   $\Gold_1^{\,1}$,
   with standard generator inclusion for one of them (say $H$) given  by $(c, x, y)_H \mapsto (c, x, d^{2}y)_G$,
   and quotient of  $G$ by the core of each subgroup being isomorphic to $S_4$.
\item Index $2$ normal subgroup $H$ of type $\Gold_2$, with standard generator 
   inclusion given by $(c, d, x, y)_H$ $\mapsto (cd, d^2, x, y^{-1})_G$, and conjugation by $c$ in
    $G$ inducing an automorphism of  $H$ that takes $(c,d,x,y)_H$ to $(cd,d,x^{-1},cy^{-1})_H$. 
\end{itemize}
\bigskip


\noindent 
\begin{tabular}{|l|}
\hline
\hline
\\[-10pt] 
$\TypeSemicolon  \Gold_2^{\,2},\>\>\>  G = \langle\, c, d, x, y \mid 
c^2, d^4, (cd)^2, x^3, [c,x], [d^2,x], d^{-1}xdx, y^3, (dy^{-1})^2, cydy \,\rangle$ \\
\\[-10pt] 
\hline
\\[-10pt] 
 $A =  \langle c,d,x \rangle \cong C_3 \times D_4$ and 
 $B = \langle c,d,y \rangle \cong S_4$, of order $24$; \ 
 $C = \langle c,d  \rangle \cong D_4$ , $\>\> (s_u,s_v) = (3,4)$
 \\[-10pt] 
${}$ \\
\hline\hline
\end{tabular}
\bigskip

\noindent
List of edge-transitive subgroups of $G$:
\medskip

\begin{itemize}[leftmargin=0.60cm]\setlength{\itemsep}{8pt}
\item Conjugacy class of four index $8$ non-normal (and non-maximal) subgroups of type   $\Gold_1$, 
   with standard generator inclusion for one of them (say $H$) given  by $(x, y)_H \mapsto (x, y)_G$, and
   quotient of  $G$ by the core of each subgroup being isomorphic to $S_4$.
\item Conjugacy class of four index $4$ non-normal maximal subgroups of type   $\Gold_1^{\,2}$, 
  with standard generator inclusion for one of them (say $H$) given by $(c, x, y)_H \mapsto (c, d^{2}y,x)_G$,
   and quotient of  $G$ by the core of each subgroup being isomorphic to $S_4$.
\item Index $2$ normal subgroup $H$ of type $\Gold_2$, with standard generator 
   inclusion given by $(c, d, x, y)_H$ $\mapsto (cd, d^2, x, y^{-1})_G$, and conjugation by $c$ in
    $G$ inducing an automorphism of  $H$ that takes $(c,d,x,y)_H$ to~$(cd,d,x,cy^{-1})_H$.
\end{itemize}
\bigskip


\noindent 
\begin{tabular}{|l|}
\hline
\hline
\\[-10pt] 
$\TypeSemicolon  \Gold_2^{\,3},\>\>\>  G = \langle\, c, d, e, x, y \mid 
c^2, d^2, e^2, [c,d], [c,e], [d,e], x^3, (cx)^2, [d,x], [e,x],$ \\ ${}$ \hskip 3.6cm $
  y^3, [c,y], y^{-1}cdye, (ye)^3  \,\rangle$ \\
\\[-10pt] 
\hline
\\[-10pt] 
 $A =  \langle c,d,e,x \rangle \cong D_6 \times C_2$ and 
 $B = \langle c,d,e,y \rangle \cong A_4 \times C_2$, of order $24$; \\ ${}$\quad  
 $C = \langle c,d,e \rangle \cong C_2 \times C_2 \times C_2$, $\>\> (s_u,s_v) = (1,2)$ 
 \\[-10pt] 
${}$ \\
\hline\hline
\end{tabular}
\bigskip

\noindent
List of edge-transitive subgroups of $G$:
\medskip

\begin{itemize}[leftmargin=0.60cm]\setlength{\itemsep}{8pt}
\item Conjugacy class of four index $8$ non-normal (and non-maximal) subgroups of type   $\Gold_1$,
   with standard generator inclusion for one of them (say $H$) given  by $(x, y)_H \mapsto (x, y)_G$, and  
   quotient of  $G$ by the core of each subgroup being isomorphic to $A_4 \times C_2$.
\item Conjugacy class of four index $4$ non-normal  maximal subgroups of type   $\Gold_1^{\,2}$,
   with standard generator inclusion for one of them (say $H$) given  by $(c, x, y)_H \mapsto (c, x,y)_G$, 
   and quotient of  $G$ by the core of each subgroup being isomorphic to $A_4$.
\item Index $2$ normal subgroup $H$ of type $\Gold_2$, with standard generator 
   inclusion given by $(c, d, x, y)_H$ $\mapsto (cd, e, x, y)_G$, and conjugation by $c$ in
    $G$ inducing an automorphism of  $H$ that takes $(c,d,x,y)_H$ to~$(c,d,x^{-1},y)_H$.
\end{itemize}
\bigskip


\noindent 
\begin{tabular}{|l|}
\hline
\hline
\\[-10pt] 
$\TypeSemicolon  \Gold_2^{\,4},\>\>\>  G = \langle\, c, d, e, x, y \mid 
c^2, d^4, e^2, (cd)^2, [c,e], [d,e], x^3, [ce,x], [d,x], (ex)^2,$ \\ ${}$ \hskip 3.6cm $
  y^3, (dy^{-1})^2, cydy, [y,e] \,\rangle$ \\
\\[-10pt] 
\hline
\\[-10pt] 
 $A =  \langle c,d,e,x \rangle \cong S_3 \times D_4$ and 
 $B = \langle c,d,e,y \rangle \cong S_4 \times C_2$, of order $48$; \\ ${}$\quad  
 $C = \langle c,d,e \rangle \cong D_4 \times C_2$, $\>\> (s_u,s_v) = (3,4)$
 \\[-10pt] 
${}$ \\
\hline\hline
\end{tabular}
\bigskip

\noindent
List of edge-transitive subgroups of $G$:
\medskip

\begin{itemize}[leftmargin=0.60cm]\setlength{\itemsep}{8pt}
\item Conjugacy class of four index $16$ non-normal subgroups of type   $\Gold_1$, 
   with standard generator inclusion for one of them (say $H$) given  by $(x, y)_H \mapsto (x, y)_G$, and  
   quotient of  $G$ by the core of each subgroup being isomorphic to $S_4 \times C_2$.
\item Conjugacy class of four index $8$ non-normal subgroups of type   $\Gold_1^{\,1}$,
   with standard generator inclusion for one of them (say $H$) given by $(c, x, y)_H \mapsto (c, x, d^{2}y)_G$, 
   and  quotient of  $G$ by the core of each subgroup being isomorphic to $S_4 \times C_2$.
\item Conjugacy class of four index $8$ non-normal (and non-maximal) subgroups of type   $\Gold_1^{\,2}$,  
   with standard generator inclusion for one of them (say $H$) given by $(c, x, y)_H \mapsto (e, x, y)_G$, 
   and  quotient of  $G$ by the core of each subgroup being isomorphic to  $S_4$. 
\item Another conjugacy class of four index $8$ non-normal subgroups of type  $\Gold_1^{\,2}$,
      with standard generator inclusion for one of them (say $H$) given by $(c, x, y)_H \mapsto (cd^{2}e, cdy, x)_G$, 
   and  quotient of  $G$ by the core of each subgroup being isomorphic to $S_4 \times C_2$.
\item Conjugacy class of four index $4$ maximal non-normal subgroups of type   $\Gold_1^{\,3}$,
  with standard generator inclusion for one of them (say $H$) given by $(c, d, x, y)_H \mapsto (e,ce,x,d^{2}y)_G$,
   and  quotient of  $G$ by the core of each subgroup being isomorphic to  $S_4$.
\item Index $4$ normal (and non-maximal) subgroup $H$ of type   $\Gold_2$, with standard generator
   inclusion given  by $(c, d, x, y)_H$ $\mapsto (cd, d^2, x, y^{-1})_G$,  and conjugation by $c$ in 
    $G$ inducing an automorphism of  $H$ that takes $(c,d,x,y)_H$ to $(cd,d,x^{-1},cy^{-1})_H$.
\item Index $2$ normal subgroup $H$ of type  $\Gold_2^{\,1}$, with standard generator
   inclusion given by $(c, d, x, y)_H$ $\mapsto (c, d, x, y)_G$, and conjugation by $e$ in  $G$
    inducing an automorphism of  $H$ that takes $(c,d,x,y)_H$ to~$(c,d,x^{-1},y)_H$.
\item Index $2$ normal subgroup $H$ of type   $\Gold_2^{\,2}$, with standard generator 
   inclusion given by $(c, d, x, y)_H$ $\mapsto (ce, de, c, y)_G$, and conjugation by $e$ in  
    $G$ inducing an automorphism of  $H$ that takes $(c,d,x,y)_H$ to~$(c,d,x^{-1},y)_H$.
\item Index $2$ normal subgroup $H$ of type   $\Gold_2^{\,3}$, with standard generator 
   inclusion given by $(c, d, e, x, y)_H$ $\mapsto (e, cde, d^2, x, y^{-1})_G$, and conjugation 
    by $c$ in  $G$ inducing an automorphism of  $H$ that takes $(c,d,e,x,y)_H$
    to $(c,de,e,x^{-1},cdy^{-1})_H$. 
\end{itemize}
\bigskip


\noindent 
\begin{tabular}{|l|}
\hline
\hline
\\[-10pt] 
$\TypeSemicolon  \Gold_3,\>\>\>  G = \langle\, c, d, x, y \mid 
c^2, d^4, (cd)^2, x^3, (dx^{-1})^2, cxdx, y^3, (dy^{-1})^2, cdydy \,\rangle$ \\
\\[-10pt] 
\hline
\\[-10pt] 
 $A =  \langle c,d,x \rangle \cong S_4$ and 
 $B = \langle c,d,y \rangle \cong S_4$, of order $24$; \ $C = \langle c,d,e \rangle \cong D_4$, $\>\> (s_u,s_v) = (4,4)$
 \\[-10pt] 
${}$ \\
\hline\hline
\end{tabular}
\medskip

\noindent 
Again here we point out that the final relator $cdydy$ corrects a typographical error 
in the description of this group in \cite{CMMP}, where it was given as $cdyd$, as noted in the Introduction.
\smallskip

\noindent
List of edge-transitive subgroups of $G$:
\medskip

\begin{itemize}[leftmargin=0.60cm]\setlength{\itemsep}{8pt}
\item Two conjugacy classes of eight index $8$ non-normal maximal subgroups isomorphic
  to  $\Gold_1$, with standard generator inclusions for a representative $H$ of each class 
   given by  $(x,y)_H \mapsto (x,y)_G$ and $(x, y)_H \mapsto (x, d^{-1}y^{-1})_G$,  and quotient of  $G$
    by the core of each subgroup 
    being isomorphic to $\PSL(2,7)$.
\end{itemize}

\smallskip
\noindent   
Note : These two conjugacy classes of subgroups (of type $\Gold_1$) fuse into a single conjugacy 
class within a group of type $\Gold_3^{\, 1}$  as well as within in a group of type $\DjM_4^{\, 2}$,
but not within a group of type $\DjM_4^{\, 1}$.

\smallskip\bigskip


\noindent 
\begin{tabular}{|l|}
\hline
\hline
\\[-10pt] 
$\TypeSemicolon  \Gold_3^{\,1},\>\>\>  G = \langle\, c, d, e, x, y \mid 
c^2, d^4, e^2, (cd)^2, [c,e], [d,e], x^3, (dx^{-1})^2, cxdx, [x,e],$ \\ ${}$ \hskip 3.6cm $
  y^3, (dy^{-1})^2, cd^{-1}ydy, [y,ed^2]  \,\rangle$ \\
\\[-10pt] 
\hline
\\[-10pt] 
 $A =  \langle c,d,x \rangle \cong S_4 \times C_2$ and 
 $B = \langle c,d,y \rangle \cong S_4 \times C_2$, of order $48$; \\ ${}$\quad  
 $C = \langle c,d,e \rangle \cong D_4 \times C_2$, $\>\> (s_u,s_v) = (5,5)$ 
 \\[-10pt] 
${}$ \\
\hline\hline
\end{tabular}
\bigskip

\noindent
List of edge-transitive subgroups of $G$:
\medskip

\begin{itemize}[leftmargin=0.60cm]\setlength{\itemsep}{8pt}
\item Conjugacy class of sixteen index $16$ non-normal subgroups isomorphic to  $\Gold_1$, 
   with standard generator inclusion for one of them (say $H$) given by $(x, y)_H \mapsto (x, y)_G$,  
   and quotient of  $G$ by the core of each subgroup being isomorphic 
   to the wreath product $\PSL(2,7) \wr  C_2$.
\item Index $2$ normal subgroup $H$ of type   $\Gold_3$, with standard generator
   inclusion given by $(c, d, x, y)_H$ $\mapsto (cd^{-1},d,y,x)_G$, and conjugation by $c$
    in  $G$ inducing an automorphism of  $H$ that takes $(c,d,x,y)_H$ 
    to $(c,d,cdx,y)_H$.
\end{itemize}

\bigskip


\noindent 
\begin{tabular}{|l|}
\hline
\hline
\\[-10pt] 
$\TypeSemicolon  \Gold_4,\>\>\>  G = \langle\, a, b, s, x, y \mid 
a^4, b^4, s^2, [a,b], sasb^{-1}, sbsa^{-1}, x^3, x^{-1}axb^{-1}, 
x^{-1}bxba, (xs)^2,$ \\ ${}$ \hskip 3.6cm $ 
  y^3, y^{-1}absysa^2, y^{-1}sa^2yba^{-1}, a^{-1}sysay  \,\rangle$ \\
\\[-10pt] 
\hline
\\[-10pt] 
 $A =  \langle a,b,s,x \rangle$ and 
 $B = \langle a,b,s,y \rangle$, of order $96$; \ $C = \langle a,b,s \rangle$, of order $32$, $\>\> (s_u,s_v) = (5,6)$
 \\[-10pt] 
${}$ \\
\hline\hline
\end{tabular}
\bigskip

\noindent
List of edge-transitive subgroups of $G$:
\medskip

\begin{itemize}[leftmargin=0.60cm]\setlength{\itemsep}{8pt}
\item Two conjugacy classes of thirty-two index $32$ non-normal maximal subgroups
   of type   $\Gold_1$, with standard generator inclusions for a representative $H$
   of each class  given  by $(x, y)_H \mapsto (x, y)_G$  
      and $(x, y)_H \mapsto (x,ay)_G$, 
   and quotient of  $G$ by the core of each subgroup being isomorphic to the 
   alternating group $A_{32}$.
\end{itemize}

\smallskip
\noindent 
Note : These two conjugacy classes of subgroups (of type $\Gold_1$) fuse into a single conjugacy 
within a group of type $\Gold_4^{\, 1}$.

\smallskip\bigskip


\noindent 
\begin{tabular}{|l|}
\hline
\hline
\\[-10pt] 
$\TypeSemicolon  \Gold_4^{\,1},\>\>\>  G = \langle\, a, b, s, t, x, y \mid 
a^4, b^4, s^2, t^2, [a,b], sasb^{-1}, sbsa^{-1}, (ta)^2, (tb)^2, [s,t],$ \\ ${}$ \hskip 1.6cm $
  x^3, x^{-1}axb^{-1}, x^{-1}bxba, (xs)^2, [x,t],$ \\ ${}$ \hskip 1.6cm $
  y^3, y ^{-1}absysa^2, y^{-1}sa^2yba^{-1}, a^{-1}sysay, tytsa^2y^{-1}a^2s  \,\rangle$ \\
\\[-10pt] 
\hline
\\[-10pt] 
 $A =  \langle a,b,s,t,x \rangle$ and 
 $B = \langle a,b,s,t,y \rangle$, of order $192$; \ $C = \langle a,b,s,t \rangle$, of order $64$, $\>\> (s_u,s_v) = (7,7)$
 \\[-10pt] 
${}$ \\
\hline\hline
\end{tabular}
\bigskip

\noindent
List of edge-transitive subgroups of $G$:
\medskip

\begin{itemize}[leftmargin=0.60cm]\setlength{\itemsep}{8pt}
\item Conjugacy class of sixty-four index $64$ non-normal subgroups isomorphic of type $\Gold_1$, 
  with standard generator inclusion for one of them (say $H$) given  by $(x, y)_H \mapsto (x, y)_G$, and 
   quotient of  $G$ by the core of each subgroup being isomorphic to the wreath   product $A_{32} \wr  C_2$.
\item Index $2$ normal subgroup $H$ of type $\Gold_4$, with standard generator \
   inclusion given by $(a, b, s, x, y)_H$ $\mapsto (a, b, s, x, y)_G$, and conjugation by $t$ 
    in  $G$ inducing an automorphism of  $H$ that takes $(a,b,s,x,y)_H$ 
    to $(a^{-1},b^{-1},s,x,sa^2ya^2s)_H$.
\end{itemize}
\bigskip


\noindent 
\begin{tabular}{|l|}
\hline
\hline
\\[-10pt] 
$\TypeSemicolon  \Gold_5,\>\>\>  G = \langle\, a, b, s, t, x, y \mid 
a^4, b^4, s^2, t^2, [a,b], sasb^{-1}, (ta)^2, (tb)^2, [s,t],$ \\ ${}$ \hskip 1.6cm $ 
  x^3, x^{-1}axb^{-1}, x^{-1}bxba, (xs)^2, [x,t],$ \\ ${}$ \hskip 1.6cm $ 
  y^3, y^{-1}absy ba^{-1}, y^{-1}sa^2ysb^{-1}a^{-1}, y^{-1}tsa^2yb^{-1}a^{-1}, y^{-1}abyb^{-1}ast, tb^{-1}ybty  \,\rangle$ \\
\\[-10pt] 
\hline
\\[-10pt] 
 $A =  \langle a,b,s,t,x \rangle$ and 
 $B = \langle a,b,s,t,y \rangle$, of order $192$; \ $C = \langle a,b,s,t \rangle$, of order $64$, $\>\> (s_u,s_v) = (6,5)$
 \\[-10pt] 
${}$ \\
\hline\hline
\end{tabular}
\bigskip

\noindent
List of edge-transitive subgroups of $G$:
\medskip

\begin{itemize}[leftmargin=0.60cm]\setlength{\itemsep}{8pt}
\item Three conjugacy classes of sixty-four index $64$ non-normal maximal subgroups 
   of type $\Gold_1$, with standard generator inclusions for a representative $H$
   of each class given  by $(x, y)_H \mapsto~(x, y)_G$, 
   $(x, y)_H \mapsto (b^{2}x,y)_G$ and $(x, y)_H \mapsto (b^{-1}x,y)_G$, and quotient of  $G$ by the core of each subgroup 
   being isomorphic to the alternating group $A_{64}$.
\item A fourth conjugacy class of sixty-four index $64$ non-normal (and non-maximal) subgroups
   of type $\Gold_1$, 
   with standard generator inclusion for one of them (say $H$) given  by $(x, y)_H \mapsto~(bx, y)_G$, 
    and quotient of  $G$ by the core of each subgroup being isomorphic to a 
    subgroup of index $3$ in the wreath product $A_4 \wr A_{16}$.
\item Conjugacy class of sixteen index $16$ non-normal maximal subgroups of type  $\Gold_2$,  
   with standard generator inclusion for one of them (say $H$) given  by $(c, d, x, y)_H \mapsto (ab^{-1}s, y^{-1}ab^{-1}sy, bx, y)_G$, 
    and quotient of  $G$ by the core of each subgroup being isomorphic to $A_{16}$.
\end{itemize}

\smallskip
\noindent
Note: The first two conjugacy classes of subgroups of type $\Gold_1$ (in the first bullet point above) 
fuse into a single conjugacy class within a group of type $\Gold_5^{\, 1}$.

\smallskip\bigskip


\noindent 
\begin{tabular}{|l|}
\hline
\hline
\\[-10pt] 
$\TypeSemicolon  \Gold_5^{\,1},\>\>\>  G = \langle\, a, b, s, t, v, x, y \mid 
a^4, b^4, s^2, t^2, v^2, [a,b], sasb^{-1}, (ta)^2, (tb)^2, [s,t],$ \\ ${}$ \hskip 1.6cm $ 
 vavb^{-2}a^{-1}, vbva^{-2}b, (vs)^2t, [v,t], x^3, x^{-1}axb^{-1}, x^{-1}bxba, (xs)^2, [x,t], [x,v],$ \\ ${}$ \hskip 1.6cm $ 
 y^3, y^{-1}absyba^{-1}, y^{-1}sa^2ysb^{-1}a^{-1}, y^{-1}tsa^2yb^{-1}a^{-1},y^{-1}abyb^{-1}ast,
   tb^{-1}ybty, $ \\ ${}$ \hskip 1.6cm $ 
 b^{-1}tvbtvtb^3a^2tb^{-1}, [vtb,y]  \,\rangle$ \\
\\[-10pt] 
\hline
\\[-10pt] 
 $A =  \langle a,b,s,t,v,x \rangle$ and 
 $B = \langle a,b,s,t,v,y \rangle$, of order $384$; \\ ${}$\quad  
 $C = \langle a,b,s,t, v\rangle$, of order $128$, $\>\> (s_u,s_v) = (8,7)$
 \\[-10pt] 
${}$ \\
\hline\hline
\end{tabular}
\bigskip

\noindent
List of edge-transitive subgroups of $G$:
\medskip

\begin{itemize}[leftmargin=0.60cm]\setlength{\itemsep}{8pt}
\item Conjugacy class of $128$ index $128$ non-normal subgroups of type   $\Gold_1$, 
   with standard generator inclusion for one of them (say $H$) given by $(x, y)_H \mapsto (x,y)_G$, and
   quotient of  $G$ by the core of each subgroup being isomorphic to the wreath 
   product $A_{64} \wr C_2$.
\item Another conjugacy class of $64$ index $128$ non-normal subgroups of type   $\Gold_1$, 
   with standard generator inclusion for one of them (say $H$) given by $(x, y)_H \mapsto (bx,y)_G$, and 
   quotient of  $G$ by the core of each subgroup being isomorphic to a 
   subgroup of index $3 \cdot 2^{15}$ in the wreath product $S_4 \wr A_{16}$.
\item A third conjugacy class of $64$ index $128$ non-normal subgroups of type  $\Gold_1$, 
   with standard generator inclusion for one of them (say $H$) given by $(x, y)_H \mapsto (b^{-1}x, y)_G$,
   and quotient of  $\Gold_5^{\,1}$ by the core of each subgroup being isomorphic to 
   $A_{64} \times C_2$.
\item Conjugacy class of $64$ index $64$ non-normal maximal subgroups of type $\Gold_1^{\,2}$, 
   with standard generator inclusion for one of them (say $H$) given by $(c,x,y)_H \mapsto (vt,ayb,b^{2}x)_G$,  
   and quotient of  $G$ by the core of each subgroup being isomorphic to $A_{64}$.
   Also each subgroup in this class is the normaliser in  $G$ of a  subgroup of index $128$ described in
   the bullet point immediately above.
\item Another conjugacy class of $64$ index $64$ non-normal subgroups of type  $\Gold_1^{\,2}$, 
   with standard generator inclusion for one of them (say $H$) given by $(c,x,y)_H \mapsto (vt, yab, a^{2}x)_G$,  
   and quotient of  $G$ by the core of each subgroup being 
   isomorphic to a subgroup of index $3 \cdot 2^{15}$ in the wreath product $S_4 \wr A_{16}$, 
   as in the second bullet point above. 
   Also each subgroup in this class is the normaliser in  $G$ of a subgroup of index $128$ described in
   the second bullet point above. 
\item Conjugacy class of $16$ index $32$ non-normal subgroups of type   $\Gold_2$, 
   with standard generator inclusion for one of them (say $H$) given by $(c,d,x,y)_H \mapsto  (st, a^{2}b^{2}t, ax, (bya)^{-1})_G$,  
   and quotient of  $G$ by the core of each subgroup  being isomorphic to  $A_{16} \times C_2$. 
\item Conjugacy class of $16$ index $16$ non-normal maximal subgroups of type   $\Gold_2^{\,2}$, 
  with standard generator inclusion for one of them (say $H$) given by $(c,d,x,y)_H \mapsto (v, a^{2}b^{2}stv, a^{2}x^{-1}, asy^{-1}a)_G$, 
  and quotient of  $G$ by the core of each subgroup being isomorphic to $A_{16}$.
\item Index $2$ normal subgroup $H$ of type  $\Gold_5$, with standard generator 
   inclusion  $(a, b, s, t, x, y)_H  \mapsto  (a, b, s, t, x, y)_G$, and conjugation by    
   $v$ in  $G$ inducing an automorphism of the subgroup $H$ that takes $(a,b,s,t,x,y)$ to 
   $(ab^2, b^{-1}a^2, st, t, x, tbyb^{-1}t)$. 
\end{itemize}

\subsection{Inclusions within Djokovi{\' c}-Miller classes}
\label{sec:InclusionsGtoDM}

${}$

\par

In this subsection, we consider each of the seven Djokovi{\' c}-Miller classes in turn,
describing each class by specifying a finite presentation for a representative $G$, 
which we again call the standard presentation of the class.

In each case, $G$ is isomorphic to the amalgamated free product $A*_C B$, where
$A=G_u$, $B=G_{\{u,v\}}$ and $C= A \cap B = G_{uv}$, for a fixed arc $(u,v)$ of $\Tree_3$.
We list generators for those subgroups $A$, $B$ and $C$, and then list a representative of 
each conjugacy class of edge-transitive subgroups of $G$, 
and provide information about the sizes of those conjugacy classes and 
their Goldschmidt or Djokovi\'c-Miller types. 
Note that if a subgroup $H\le G$ has Goldschmidt or Djokovi\'c-Miller type $\cH$, 
 then it can be presented in terms of the standard presentation for the class $\cH$, 
and then the inclusion of $H$ in $G$ is described by specifying the images of the standard 
generators for $\cH$ as words in the standard generators for $\cG$. 
The standard inclusions of these subgroups are given in the same manner as in the previous subsection.

To find out which of the Goldschmidt classes are included in which of the Djokovi{\' c}-Miller classes,
again we made extensive use of {\sc Magma}, in the obvious way, analogous to the approach taken in the previous subsection (but without the need for special treatment when the order of $G_{uv}$ is large).  
\medskip


\noindent 
\begin{tabular}{|l|}
\hline
\hline
\\[-10pt] 
$\TypeSemicolon  \DjM_1,\>\>\>  G = \langle\, h, a \mid h^3, a^2 \,\rangle$     \\
\\[-10pt] 
\hline
\\[-10pt] 
 $A = \langle h \rangle \cong C_3$ of order $3$; \ 
 $B = \langle a \rangle \cong C_2$ of order $2$; \ 
 $C = \{1\}$ (trivial) 
 \\[-10pt] 
${}$ \\
\hline\hline
\end{tabular}
\bigskip

\noindent
List of edge-transitive  subgroups of $G$:
\medskip

\begin{itemize}[leftmargin=0.60cm]\setlength{\itemsep}{8pt}
\item Index $2$ semisymmetric normal subgroup $H$ of type   $\Gold_1$, with standard generator 
  inclusion given by $(x, y)_H  \mapsto (h,aha)_G$,  and conjugation by $a$ in $G$
 inducing an automorphism of  $H$ that takes $(x,y)_H$ to $(y,x)_H$.
\end{itemize}
\bigskip


\noindent 
\begin{tabular}{|l|}
\hline
\hline
\\[-10pt] 
$\TypeSemicolon  \DjM_2^{\,1},\>\>\>  G = \langle\, h, p, a \mid h^3, p^2, a^2, (hp)^2, [a,p] \,\rangle$     \\
\\[-10pt] 
\hline
\\[-10pt] 
 $A = \langle h, p \rangle \cong S_3$ of order $6$; \ 
 $B = \langle p, a \rangle \cong C_2 \times C_2$ of order $4$; \ 
 $C = \langle p \rangle \cong C_2$ 
 \\[-10pt] 
${}$ \\
\hline\hline
\end{tabular}
\bigskip

\noindent
List of edge-transitive  subgroups of $G$:
\medskip

\begin{itemize}[leftmargin=0.60cm]\setlength{\itemsep}{8pt}
\item Index $4$ semisymmetric normal (and non-maximal) subgroup $H$ of type   $\Gold_1$, with standard generator
  inclusion given by $(x, y)_H \mapsto (h,aha)_G$,  and conjugation by $a$ and $p$ in $G$
  inducing automorphisms of  $H$ that take $(x,y)_H$ to $(y,x)_H$ and $(x^{-1},y^{-1})_H$.
\item Index $2$ semisymmetric normal subgroup $H$ of type   $\Gold_1^{\,1}$, with standard generator 
   inclusion given by $(c, x, y)_H \mapsto (p,h,aha)_G$,  and conjugation by $a$ in $G$  
   inducing an automorphism of  $H$ that takes $(c,x,y)_H$ to $(c,y,x)_H$.
\item Two index $2$  arc-transitive normal subgroups  of type $\DjM_1$, with standard   
   generator inclusions given by $(h,a)_H \mapsto (h,a)_G$ and $(h,a)_H \mapsto (h,ap)_G$,   
   and conjugation by $p$ in $G$  inducing an  automorphism of $H$ that  
   takes $(h,a)_H$ to $(h^{-1},a)_H$  in both cases.
\end{itemize}

\smallskip
\noindent
Note:  The last two subgroups (of type $\DjM_1$) are conjugate within a group of type $\DjM_3$.

\smallskip\bigskip


\noindent
\begin{tabular}{|l|}
\hline
\hline
\\[-10pt] 
$\TypeSemicolon  \DjM_2^{\,2},\>\>\>  G = \langle\, h, p, a \mid h^3, p^2, pa^{2}, (hp)^2 \,\rangle$     \\
\\[-10pt] 
\hline
\\[-10pt] 
 $A = \langle h, p \rangle \cong S_3$ of order $6$; \ 
 $B = \langle a \rangle \cong C_4$; \ 
 $C = \langle p \rangle \cong C_2$ 
  \\[-10pt] 
${}$ \\
\hline\hline
\end{tabular}
\bigskip

\noindent
List of edge-transitive  subgroups of $G$:
\medskip

\begin{itemize}[leftmargin=0.60cm]\setlength{\itemsep}{8pt}
\item Index $4$ semisymmetric normal (and non-maximal) subgroup $H$ of type   $\Gold_1$, with standard generator    
  inclusion given by $(x, y)_H \mapsto (h,aha^{-1})_G$,  and conjugation by $a^{-1}$ and $p$    
  in $G$ inducing automorphisms of  $H$ that take $(x,y)_H$ to $(y,x^{-1})_H$   
  and $(x^{-1},y^{-1})_H$ respectively. 
\item Index $2$ semisymmetric normal subgroup $H$ of type   $\Gold_1^{\,1}$, with standard generator    
   inclusion given by $(c, x, y)_H \mapsto (p,h,aha^{-1})_G$,  and conjugation by $a^{-1}$ in   
   $G$ inducing an automorphism of  $H$ that takes $(c,x,y)_H$ to $(c,y,x^{-1})_H$. 
\end{itemize}
\bigskip


\noindent
\begin{tabular}{|l|}
\hline
\hline
\\[-10pt] 
$\TypeSemicolon  \DjM_3,\>\>\>  G = \langle\, h, p, q, a \mid h^3, p^2, q^2, a^2, [p,q], apaq, [h,p], (hq)^2 \,\rangle$     \\
\\[-10pt] 
\hline
\\[-10pt] 
 $A = \langle h, p, q \rangle \cong S_3 \times C_2$ of order $12$; \ 
 $B = \langle a, p, q \rangle \cong D_4$ of order $8$; \\
 $C = \langle p,q \rangle \cong C_2 \times C_2 $ of order $4$ 
 \\[-10pt] 
${}$ \\
\hline\hline
\end{tabular}
\bigskip

\noindent
List of edge-transitive  subgroups of $G$:
\medskip

\begin{itemize}[leftmargin=0.60cm]\setlength{\itemsep}{8pt}
\item Index $8$ semisymmetric normal (and non-maximal) subgroup $H$ of type   $\Gold_1$, with standard generator   
   inclusion given by $(x, y)_H \mapsto (h,aha)_G$,  and conjugation by $a$, $p$ and $q$ in   
   $G$ inducing automorphisms of  $H$ that take $(x,y)_H$ to $(y,x)_H$, $(x,y^{-1})_H$   
   and $(x^{-1},y)_H$, respectively.
\item Index $4$ semisymmetric normal (and non-maximal) subgroup $H$ of type   $\Gold_1^{\,1}$, with standard generator  
  inclusion given by $(c, x, y)_H \mapsto (pq,h,aha)_G$,  and conjugation by $a$, $p$ and~$q$   
  in $G$ inducing automorphisms of  $H$ that take $(c,x,y)_H$ to $(c,y,x)_H$,  
   $(c,x,y^{-1})_H$ and $(c,x^{-1},y)_H$, respectively.
\item Conjugacy class of two index $4$ semisymmetric subgroups of type   $\Gold_1^{\,2}$,  
   with standard generator inclusion for one of them (say $H$) given by   
    $(c, x, y)_H \mapsto (p,aha,h)_G$ and the quotient of $G$ by the core of each subgroup being isomorphic to $D_4$.
\item Index $2$ semisymmetric normal subgroup $H$ of type   $\Gold_1^{\,3}$,  with standard generator  
   inclusion given by $(c, d, x, y)_H \mapsto (q,p,h,aha)_G$,  and conjugation by $a$ in  
    $G$ inducing an automorphism of  $H$ that takes $(c,d,x,y)_H$ to $(d,c,y,x)_H$.
\item Conjugacy class of two index $4$ arc-transitive non-normal subgroups of type   $\DjM_1$,  
   with standard generator inclusion for one of them (say $H$) given by   
    $(h,a)_H \mapsto (h,a)_G$,  and     
   quotient of $G$ by the core of each subgroup being isomorphic to $D_4$.
\item Index $2$ arc-transitive normal subgroup $H$ of type   $\DjM_2^{\,1}$,  with standard generator  
   inclusion given by $(h,p,a)_H$ $\mapsto (h,pq,a)_G$,  and conjugation by $p$ in~$G$   
    inducing an automorphism of $H $ that takes $(h,p,a)_H$ to $(h,p,ap)_H$.
\item Index $2$ arc-transitive normal subgroup $H$ of type   $\DjM_2^{\,2}$,  with standard generator  
   inclusion given by $(h,p,a)_H$ $\mapsto (h,pq,ap)_G$,  and conjugation by $p$ in~$G$   
    inducing an automorphism of $H$ that takes $(h,p,a)_H$ to $(h,p,ap)_H$.
\end{itemize}
\bigskip


\noindent
\begin{tabular}{|l|}
\hline
\hline
\\[-10pt] 
$\TypeSemicolon  \DjM_4^{\,1},\>\>\>  G = \langle\, h, p, q, r, a \mid h^3, p^2, q^2, r^2, a^2, 
[p,q], [p,r], p(qr)^2, [a,p], aqar, h^{-1}phq, $ \\ ${}$ \qquad\qquad\qquad $h^{-1}qhpq, (hr)^2 \,\rangle$     \\
\\[-10pt] 
\hline
\\[-10pt] 
 $A = \langle h, p, q, r \rangle \cong S_4$ of order $24$; \ 
 $B = \langle a, p, q, r \rangle \cong D_4 \rtimes C_2$ of order $16$; \\
 $C = \langle p,q,r \rangle \cong D_4 $ of order $8$ 
 \\[-10pt] 
${}$ \\
\hline\hline
\end{tabular}
\bigskip

\noindent
List of edge-transitive  subgroups of $G$:
\medskip

\begin{itemize}[leftmargin=0.60cm]\setlength{\itemsep}{8pt}
\item Two conjugacy classes of eight index $16$ semisymmetric non-normal subgroups of type $\Gold_1$,  
   with standard generator inclusion for a representative $H$ of each class given by   
    $(x, y)_H \mapsto (h,aha)_G$ and  $(x, y)_H \mapsto (ph,aha)_G$,
     and the quotient of $G$ by the core of each subgroup being isomorphic to $\PGL(2,7)$.
\item Index $2$ semisymmetric normal subgroup $H$ of type   $\Gold_3$, with standard generator    
   inclusion given by $(c, d, x, y)_H$ $\mapsto (pr,pqr,h^{-1},aha)_G$,  and conjugation  
   by $a$ in $G$ inducing an automorphism of  $H$
   that takes $(c,d,x,y)_H$    
   to $(cd,d^{-1},y^{-1},x^{-1})_H$.
\item Two conjugacy classes of eight index $8$ arc-transitive non-normal maximal subgroups of type $\DjM_1$, 
    with standard generator inclusion for one of them (say $H$) 
   given by  $(h,a)_H \mapsto (h,a)_G$ and $(h,a)_H \mapsto (h,ap)_G$,  and quotient of $G$   
    by the core of each subgroup (in each class) being isomorphic to $\PGL(2,7)$.
\end{itemize}

\smallskip 
\noindent
Note: The two conjugacy classes of subgroups of type $\Gold_1$ fuse into a single 
conjugacy class within a group of type $\DjM_5$, as do the two conjugacy classes 
of subgroups of type $\DjM_1$.

\smallskip\bigskip


\noindent
\begin{tabular}{|l|}
\hline
\hline
\\[-10pt] 
$\TypeSemicolon  \DjM_4^{\,2},\>\>\>  G = \langle\, h, p, q, r, a \mid h^3, p^2, q^2, r^2, a^4, a^{2}p, 
[p,q], [p,r], p(qr)^2, a^{-1}qar, h^{-1}phq, $ \\ ${}$ \qquad\qquad\qquad $h^{-1}qhpq, (hr)^2 \,\rangle$     \\
\\[-10pt] 
\hline
\\[-10pt] 
 $A = \langle h, p, q, r \rangle \cong S_4$ of order $24$; \ 
 $B = \langle a, p, q, r \rangle \cong C_8 \rtimes_3 C_2$ of order $16$; \\
 $C = \langle p,q,r \rangle \cong D_4 $ of order $8$ 
 \\[-10pt] 
${}$ \\
\hline\hline
\end{tabular}
\bigskip

\noindent
List of edge-transitive  subgroups of $G$:
\medskip

\begin{itemize}[leftmargin=0.60cm]\setlength{\itemsep}{8pt}
\item Conjugacy class of $16$ index $16$ semisymmetric subgroups of type   $\Gold_1$,  
  with standard generator inclusion for one of them (say $H$) given by  
  $(x, y)_H \mapsto (aha^{-1},h)_G$, and quotient of $G$ by the core of each subgroup being isomorphic to the
  wreath product $\PSL(2,7)\wr C_2$.
\item Index $2$ semisymmetric normal subgroup $H$ of type $\Gold_3$, with standard generator    
  inclusion  given by $(c, d, x, y)_H$ $\mapsto (pr,pqr,h^{-1},aha^{-1})_G$,  and conjugation   
  by $a^{-1}$ in $G$ inducing an automorphism of   $H$
   that takes $(c,d,x,y)_H$   
  to $(cd,d^{-1},y^{-1},x^{-1})_H$. 
\end{itemize}
\bigskip


\noindent
\begin{tabular}{|l|}
\hline
\hline
\\[-10pt] 
$\TypeSemicolon  \DjM_5,\>\>\>  G = \langle\, h, p, q, r, s, a \mid h^3, p^2, q^2, r^2, s^2, a^2, 
[p,q], [p,r], [p,s], [q,r], [q,s], pq(rs)^2,  $ \\ ${}$ \qquad\qquad\qquad $apaq, aras, [h,p], h^{-1}qhr, h^{-1}rhpqr, (hs)^2 \,\rangle$     \\
\\[-10pt] 
\hline
\\[-10pt] 
 $A = \langle h, p, q, r, s \rangle \cong S_4 \times C_2$ of order $24$; \ 
 $B = \langle a, p, q, r, s \rangle \cong C_8 \rtimes V_4$ of order $32$; \\
 $C = \langle p,q,r,s \rangle \cong D_4 \times C_2$ of order $16$ 
 \\[-10pt] 
${}$ \\
\hline\hline
\end{tabular}
\bigskip

\noindent
List of edge-transitive  subgroups of $G$:
\medskip

\begin{itemize}[leftmargin=0.60cm]\setlength{\itemsep}{8pt}
\item Conjugacy class of sixteen index $32$ semisymmetric non-normal subgroups of type   $\Gold_1$,   
   with standard generator inclusion for one of them (say $H$) given by  $(x, y)_H \mapsto (h,aha)_G$. 
  The quotient of $G$ by the core of $H$ is isomorphic to the automorphism group of the Biggs-Conway 
  graph, namely a $C_2$-extension of the wreath product $\PSL(2,7) \wr C_2$,  
  isomorphic to a semi-direct product of $\PSL(2,7) \times \PSL(2,7)$ by $C_2 \times C_2$. 
\item Index $4$ semisymmetric normal (and non-maximal) subgroup $H$ of type  $\Gold_3$, with standard generator  
  inclusion  given by $(c, d, x, y)_H$ $\mapsto (pr,sr,aha,h^{-1})_G$,  and conjugation   
  by $a$ and $p$ in $G$ inducing automorphisms of  $H$ that take $(c,d,x,y)_H$   
  to $(cd,d^{-1},y^{-1},x^{-1})_H$ and $(c,d,cdx,y)_H$, respectively.
\item Index $2$ semisymmetric normal subgroup $H$ of type  $\Gold_3^{\,1}$, with standard generator  
  inclusion  given by $(c, d, e, x, y)_H$ $\mapsto (qs, sr, p, h^{-1}, aha)_G$,  and conjugation   
  by $a$ in~$G$ inducing an automorphism of  $H$ that takes $(c,d,e,x,y)_H$   
  to $(cd,d^{-1},d^{2}e,y^{-1},x^{-1})_H$. 
\item Conjugacy class of sixteen index $16$ arc-transitive non-normal subgroups of type
  $\DjM_1$, with standard generator inclusion for one of them (say $H$) given by   
    $(h,a)_H \mapsto (h,a)_G$,  and quotient of $G$ by the core of each subgroup   
    being isomorphic to a subgroup of index $2$ in $\PGL(2,7) \wr C_2$.
\item Index $2$ arc-transitive normal subgroup $H$ of type  $\DjM_4^{\,1}$, with standard generator   
   for one of them (say $H$) given by $(h,p,q,r,a)_H \mapsto (h^{-1},pq,qr,ps,a)_G$,  
   and conjugation  by $p$ in $G$ inducing an automorphism of $H$ that takes $(h,p,q,r,a)_H$   
   to $(h,p,q,r,ap)_H$.
\item Index $2$ arc-transitive normal subgroup $H$ of type  $\DjM_4^{\,2}$, with standard generator   
   inclusion given by $(h,p,q,r,a)_H$ $\mapsto (h^{-1},pq,qr,ps,ap)_G$,  and conjugation  
   by $p$ in $G$ inducing an automorphism of $H$ that takes $(h,p,q,r,a)_H$  
   to $(h,p,q,r,ap)_H$.
\end{itemize}

\section{Finding all edge-transitive cubic graphs of order up to 10000}
\label{sec:10000list}

All connected arc-transitive cubic graphs on up to $10000$ vertices were 
determined by the first author in 2011 (see \cite{ATConder}),
 but this list was not widely publicised before now.

The graphs were found with the help of an improved version of the {\tt LowIndexNormalSubgroups} command 
in {\sc Magma}, using the same approach as was used to find all of order up to $2048$ 
five years earlier.  For example, those of type $\DjM_1$ were determined by \\[+2pt]
(a) using {\tt LowIndexNormalSubgroups} to find all normal subgroups $N$ of index 
up to $30000$ in the corresponding finitely-presented group $G$, such that the orders of the 
images in $G/N$ of the stabilisers $A=G_u$, $B=G_{\{u,v\}}$ and $C= A \cap B = G_{uv}$ are preserved, 
and then for each such $N$, \\[+2pt]
(b) constructing the associated arc-transitive cubic graph $\Gamma$ on which $G/N$ has an arc-transitive 
action of type $\DjM_1$, and \\[+2pt]
(c) checking that the full automorphism group of $\Gamma$ has order $|G/N|$. 

Step (c) can be performed either by using the {\tt AutomorphismGroup} command for graphs 
in {\sc Magma}, or by testing whether or not the subgroup $N$ is normal in 
one of the finitely-presented groups associated with an amalgam that contains the given one. 

The analogous process was used for each of the other six Djokovi{\' c}-Miller classes, 
 applying the {\tt LowIndexNormalSubgroups} command to find normal subgroups of index 
 up to $10000c$, where $c = 6$ for $\DjM_2^{\,1}$ and $\DjM_2^{\,2}$, or $c = 12$ for $\DjM_3$, 
 or $c = 24$ for $\DjM_4^{\,1}$ and $\DjM_4^{\,2}$, or $c = 48$ for $\DjM_5$. 

\smallskip

The total number of connected arc-transitive cubic graphs on up to $10000$ vertices is $3815$, 
with the numbers of each Djokovi{\' c}-Miller type given in the next section. 

\medskip

We subsequently used an approach analogous to the one above, in order to find all 
connected semisymmetric cubic graphs on up to $10000$ vertices.
Here we applied  {\tt LowIndexNormalSubgroups}  
to find all normal subgroups of the appropriate index in each of the finitely-presented groups
associated with the thirteen `smallest' Goldschmidt amalgams, 
namely all except $\Gold_4^{\,1}$, $\Gold_5$ and $\Gold_5^{\,1}$.

The graph construction in step (b) of the semisymmetric analogue of the above process 
differed a little from the arc-transitive version. 
First, let $\overline{G}$ denote the quotient $G/N$ of the relevant finitely-presented group $G$ 
by the normal subgroup $N$, and let $\overline{G_u} = G_{u}N/N$ and $\overline{G_v} = G_{v}N/N$ 
(which have the same orders as $G_u$ and $G_v$ and as each other).  
Then the edge-transitive graph $\Gamma$ we construct is the bipartite graph with vertex set the union 
of the two right coset spaces $(\overline{G}:\overline{G_u})$ and $(\overline{G}:\overline{G_v})$, 
and the edges are all pairs $\{\overline{G_u}\,\overline{g},\overline{G_v}\,\overline{g}\}$ 
obtained by right-multiplying the pair $\{\overline{G_u},\overline{G_v}\}$ by elements of $\overline{G}$.
Note that this graph $\Ga$ has order $2\,|\overline{G}\!:\!\overline{G_u}|$, 
and so $|G/N| = |\overline{G}| = |\overline{G_u}||V(\Ga)|/2 = |G_u||V(\Ga)|/2$. 

Otherwise steps (a) and (c) were straightforward analogues of the corresponding steps in the 
arc-transitive case. 

\smallskip

For the remaining three amalgams, namely $\Gold_4^{\,1}$, $\Gold_5$ and $\Gold_5^{\,1}$, 
the above process was challenging because the upper bound on the index $|G\!:\!N|$ 
(of $960000$, $960000$ and $1920000$ respectively) was too large.  
(The current implementation of {\tt LowIndexNormalSubgroup} command 
in {\sc Magma} computes normal subgroups of index up to at most $500000$.)
Instead we took a different approach, after observing that the finitely-presented groups associated 
with the $\Gold_4$ and $\Gold_5$ amalgams are perfect.

In the case of $\Gold_4^{\,1}$, we first consider the finitely-presented group $H$ associated 
with $\Gold_4$.  This group is perfect, and by the {\tt LowIndexNormalSubgroups} computation 
undertaken for it, $H$ has just two proper normal subgroups
of index up to $960000/2 = 480000$, namely one of index $6048$ and one of index $372000$, 
with simple quotients $\PSU(3,3)$ and $\PSL(3,5)$.  These give rise to semisymmetric graphs 
of type $\Gold_4^{\,1}$, with orders $2 \cdot 6048/96 = 126$ and $2 \cdot 372000/96 = 7750$. 

On the other hand, the {\tt SimpleQuotientProcess} in {\sc Magma} shows that the 
finitely-presented group $G$ associated with $\Gold_4^{\,1}$ has no non-abelian simple quotient 
of order up to $960000$, and then since this group $G$ has a unique subgroup of index $2$, 
isomorphic to the group $H$ considered above for $\Gold_4$, it follows that $G$ has at most 
two proper normal subgroups 
of index up to $960000$, namely one with index $12096$ and one with index $744000$.  
Indeed the {\tt Homomorphisms} 
command in {\sc Magma} shows that there is a unique homomorphism from $G$ onto 
$\Aut(\PSU(3,3)) \cong {\rm P{\Sigma}U}(3,3)$, with kernel of index $12096$, and similarly a unique 
homomorphism from $G$ onto $\Aut(\PSL(3,5)) \cong \PSL(3,5)\rtimes C_2$, with kernel 
of index $744000$. 
Hence there are exactly two edge-transitive graphs of type $\Gold_4^{\,1}$ with order at most $10000$, 
namely the two semisymmetric graphs of orders $126$ and $7750$ mentioned in the paragraph above. 

For $\Gold_5$, every quotient of the associated group $G$ up to the maximum conceivable 
order $960000$ must have a non-abelian simple quotient, and in fact an application of 
the {\tt SimpleQuotientProcess} in {\sc Magma} \cite{Magma} shows that this group $G$ has only one 
non-abelian simple quotient of order up to $960000$, namely the Mathieu group  M$_{12}$, of order $95040$. 
It follows that every non-trivial quotient $Q$ of $G$ of order at most $960000$ has order $95040m$ 
for some $m \le 10$, and hence is isomorphic to an extension by M$_{12}$ of a soluble normal subgroup 
of order at most $10$. But in that case, the normal subgroup must be central in the perfect quotient $Q$, 
and hence has order at most $2$, because the Schur multiplier of M$_{12}$ is cyclic of order $2$.
By a {\sc Magma} computation using the {\tt Darstellungsgruppe} function and the {\tt Homomorphisms} command, however, it can be shown that the Schur cover of M$_{12}$ (of order $2 \cdot 95040 = 190080$) 
is not a quotient of $G$, and hence the only possibility is that $Q \cong {\rm M}_{12}$.  

Next, the {\tt Homomorphisms} command in {\sc Magma} shows there are just two homomorphisms 
from $G$ onto M$_{12}$, up to equivalence under conjugation in $\Aut({\rm M}_{12})$, 
but these have the same kernel, so both give rise to the same edge-transitive cubic graph 
of order $2 \cdot 95040/192 = 990$, namely a semisymmetric one of type $\Gold_5^{\,1}$.  
Hence there is no edge-transitive graph of type $\Gold_5$ with order at most $10000$. 

For $\Gold_5^{\,1}$, the {\tt SimpleQuotientProcess} in {\sc Magma} shows that the associated 
finitely-presented group $G$ has no non-abelian simple quotient of order up to $1920000$, 
and then since $G$ has a unique subgroup of index $2$, isomorphic to the one considered above 
for $\Gold_5$, it follows that there is just one edge-transitive graph of type $\Gold_5^{\,1}$ 
with order at most $10000$, namely the semisymmetric graph of order $990$ 
mentioned in the paragraph above. 

\smallskip

The total number of connected semisymmetric cubic graphs on up to $10000$ vertices is $1043$, 
with the numbers of each Goldschmidt type given in the next section. 

\section{Examples of edge-transitive cubic graphs of each type}
\label{sec:Examples}

Here we give some information about examples of connected finite edge-transitive 
cubic graphs of each type. 
For most of the types, we can even give the smallest possible example. 
When mentioned, the `action type' indicates the type of action for each conjugacy class 
of arc-transitive subgroups of the automorphism group in the arc-transitive case, 
or for each conjugacy class of semisymmetric subgroups of the automorphism group in the semisymmetric case. 

\subsection{Examples with given arc-transitive type}

${}$ 

\medskip

\noindent {\bf (a) Type  $\DjM_1$}:  \\[+4pt]
The smallest example is the graph F026 in \cite{Foster,CD}.
Other small examples have orders $38$, $42$, $56$, $62$, $74$, $78$, $86$ and $98$.
There are $2522$ examples of type $\DjM_1$ with order up to $10000$.
\medskip

\noindent {\bf (b) Type  $\DjM_2^{\,1}$}:  \\[+4pt]
The smallest example is the complete graph $K_4$, and the next smallest is the $3$-cube $Q_3$,  
known as F004 and F008 in \cite{Foster,CD}, with both having action type $(\DjM_1,\DjM_2^{\,1})$. 
Other small examples have orders $16$, $24$, $32$, $48$, $50$, $54$, $60$, $64$, $72$, $84$, $96$ and $98$.
There are $1135$ examples of type $\DjM_2^{\,1}$ with order up to $10000$.
\medskip

\noindent {\bf (c) Type  $\DjM_2^{\,2}$}:  \\[+4pt]
The smallest example is the graph F448C in \cite{CD} but not in \cite{Foster}, 
with action type $(\DjM_2^{\,2})$.  The next smallest has order $896$. 
There are $9$ examples of type $\DjM_2^{\,2}$ with order up to $10000$.
\medskip

\noindent {\bf (d) Type  $\DjM_3$}:  \\[+4pt]
The smallest example is the complete bipartite graph $K_{3,3}$, and the next smallest is 
the Petersen graph, known as F006 and F010 in \cite{Foster,CD}, 
and these have action types $(\DjM_1,\DjM_2^{\,1},\DjM_2^{\,2},\DjM_3)$ and $(\DjM_2^{\,1},\DjM_3)$ 
respectively. 
Other small examples have orders $18$, $20$, $28$, $40$, $56$, $80$ and $96$.
There are 129 examples of type  $\DjM_3$ with order up to $10000$.
\medskip

\noindent {\bf (e) Type  $\DjM_4^{\,1}$}:  \\[+4pt]
The smallest example is the Heawood graph F014 in \cite{Foster,CD}, 
with action type $(\DjM_1,\DjM_4^{\,1})$.  The next smallest example is F102 in \cite{Foster,CD}.   
There are 13 examples of type $\DjM_4^{\,1}$ with order up to $10000$.
\medskip

\noindent {\bf (f) Type  $\DjM_4^{\,2}$}:  \\[+4pt]
There are no examples of type $\DjM_4^{\,2}$ with order up to $10000$.
The smallest example was shown in \cite{Cond} to be a unique one of order $5314410$
that happens to be a $3^{11}$-fold cover of Tutte's 8-cage, with action type necessarily $(\DjM_4^{\,2})$. 
A larger example (necessarily with the same action type) is a $3^{10}$-fold cover of F468, 
with order $27634932$, as mentioned in \cite{CN}. 

\medskip

\noindent {\bf (g) Type  $\DjM_5$}:  \\[+4pt]
The smallest example is Tutte's 8-cage F030 in \cite{Foster,CD}, 
and the next smallest is a triple cover of that, known as F090 in \cite{Foster,CD},   
with both having action type $(\DjM_4^{\,1},\DjM_4^{\,1},\DjM_5)$.    
There are 7 examples of type $\DjM_5$ with order up to $10000$.

\subsection{Examples with given semisymmetric  type}

${}$ 

\medskip

\noindent {\bf (a) Type  $\Gold_1$}:  \\[+4pt]
The smallest example is the Foster-Ljubljana graph of order 112, attributed to Foster 
and described in considerable detail in \cite{CMMPP}.
Other small examples were found in \cite{CMMP}, with orders 336, 378, 400, and so on.  
There are 643 examples of type  $\Gold_1$ with order up to $10000$.
\medskip

\noindent {\bf (b) Type  $\Gold_1^{\,1}$}:  \\[+4pt]
The smallest example has order 144, with action type $(\Gold_1, \Gold_1^{\,1})$. 
Other small examples found in \cite{CMMP} have orders 216, 336, 432, and so on. 
There are 222 examples of type  $\Gold_1^{\,1}$ with order up to $10000$.
\medskip

\noindent {\bf (c) Type  $\Gold_1^{\,2}$}:  \\[+4pt]
The smallest example has order 294, with action type $(\Gold_1, \Gold_1^{\,2})$. 
Another small example found in \cite{CMMP} has order 504.
There are 51 examples of type  $\Gold_1^{\,2}$ with order up to $10000$.
\medskip

\noindent {\bf (d) Type  $\Gold_1^{\,3}$}:  \\[+4pt]
The smallest example has order 120, with action type $(\Gold_1, \Gold_1^{\,1}, \Gold_1^{\,2}, \Gold_1^{\,2}, \Gold_1^{\,3})$, 
indicating that the automorphism group has two different (normal) subgroups of index 2 acting 
with type  $\Gold_1^{\,2}$. 
Other small examples found in \cite{CMMP} have orders 220, 240, 336, and so on. 
There are 90 examples of type  $\Gold_1^{\,3}$ with order up to $10000$.
\medskip

\noindent {\bf (e) Type  $\Gold_2$}:  \\[+4pt]
There are no examples of type  $\Gold_2$ with order up to $10000$.
With the help of {\sc Magma} \cite{Magma}, the smallest example can be shown to have order 25272, 
with automorphism group isomorphic to a split extension of $C_3^{\,3}$ by $\PSL(3,3)$, 
and with action type $(\Gold_2)$.  
\medskip

\noindent {\bf (f) Type  $\Gold_2^{\,1}$}:  \\[+4pt]
The smallest example has order 110, with action type $(\Gold_2, \Gold_2^{\,1})$.  
Other small examples found in \cite{CMMP} have orders 182, 330, 506 and 546. 
There are 17 examples of type  $\Gold_2^{\,1}$ with order up to $10000$.
\medskip

\noindent {\bf (g) Type  $\Gold_2^{\,2}$}:  \\[+4pt]
There are no examples of type  $\Gold_2^{\,2}$ with order up to $10000$.
The smallest example we have found has order 2527220401920, with automorphism group 
isomorphic to the simple Mathieu group M$_{24}$, and with action type $(\Gold_2^{\,2})$. 
\medskip

\noindent {\bf (h) Type  $\Gold_2^{\,3}$}:  \\[+4pt]
There are no examples of type  $\Gold_2^{\,3}$ with order up to $10000$.
We found an example with order $19!/12$, 
coming from an action of $\Sym(19)$ with type  $\Gold_2^{\,3}$ that is not extendable to an  action of a 
larger group with type  $\Gold_2^{\,4}$. The action type for this graph is $(\Gold_2, \Gold_2^{\,3})$.
The smallest example we have found has order 39366, with (soluble) automorphism group 
of order $472392 = 2^{3}\cdot 3^{10}$, and action type $(\Gold_1,\Gold_1^{\,2},\Gold_2,\Gold_2^{\,3})$. 
\medskip

\noindent {\bf (i) Type  $\Gold_2^{\,4}$}:  \\[+4pt]
The smallest example of type  $\Gold_2^{\,4}$ is the well-known Gray graph, with order 54, 
with (soluble) automorphism group of order 1296, and action type 
$(\Gold_1, \Gold_1^{\,1}, \Gold_1^{\,2}, \Gold_1^{\,2}, \Gold_1^{\,3}, \Gold_2, \Gold_2^{\,1}, \Gold_2^{\,2}, \Gold_2^{\,3}, \Gold_2^{\,4})$,
indicating that the automorphism group has many different subgroups that act semisymmetrically 
on the graph, including two different subgroups of index 8 acting with type  $\Gold_1^{\,2}$. 
There are 16 examples of type  $\Gold_2^{\,4}$ with order up to $10000$.
\medskip

\noindent {\bf (j) Type  $\Gold_3$}:  \\[+4pt]
There are no examples of type  $\Gold_3$ with order up to $10000$.
We found an easy example with order $15!/24$, 
coming from an action of $\Alt(15)$ with type  $\Gold_3$ that is not extendable to an  action of a 
larger group with type  $\Gold_3^{\,1}$, $D_4^{\,1}$ or $D_4^{\,2}$. The action type for this graph is $(\Gold_3)$.
A smaller example is a graph of order 501645312, coming from the action of a group of order 6019743744 
that is isomorphic to an extension of a group of order $2^{14}\,3^{7}$ by $\PSL(2,7)$, having a transitive 
but imprimitive permutation representation on 28 points with 7 blocks of size 4 (making it a subgroup of 
the wreath product $\Sym(4) \wr \PSL(3,2)$, but it is not isomorphic to $\Alt(4) \wr \PSL(3,2)$).  
Again this action not extendable to one of a larger group of type  $\Gold_3^{\,1}$, $D_4^{\,1}$ or $D_4^{\,2}$. 
The action type for this graph is $(\Gold_1, \Gold_3)$. 
\medskip

\noindent {\bf (k) Type  $\Gold_3^{\,1}$}:  \\[+4pt]
The smallest example of type  $\Gold_3^{\,1}$ is a graph of order 5760, 
with an (insoluble) automorphism group of order 138240, and with action type $(\Gold_3, \Gold_3^{\,1})$. 
This is the only example with order up to $10000$, 
\medskip

\noindent {\bf (l) Type  $\Gold_4$}:  \\[+4pt]
There are no examples of type  $\Gold_4$ with order up to $10000$.
We easily found an example with order $28!/96$, coming from an action of $\Alt(28)$ with type  $\Gold_4$ 
that is not extendable to an action of a larger group with type  $\Gold_4^{\,1}$. The action type for this graph is $(\Gold_4)$.
\medskip

\noindent {\bf (m) Type  $\Gold_4^{\,1}$}:  \\[+4pt]
The smallest example of type  $\Gold_4^{\,1}$ is a graph of order 126, 
with automorphism group $\Aut(\PSU(3,3))$ of order 12096, and action type $(\Gold_4, \Gold_4^{\,1})$.  
There is another of order 7750 with automorphism group $\PGL(3,5)$ of order $372000$,  
and the same action type.   
These are the only examples of type  $\Gold_4^{\,1}$ with order up to 10000. 
\vskip 5pt

\noindent {\bf (n) Type  $\Gold_5$}:  \\[+4pt]
There are no examples of type  $\Gold_5$ with order up to $10000$.
We found an easy example with order $62!/192$, 
coming from an action of $\Alt(62)$ with type  $\Gold_5$ that is not extendable to an  
action of a larger group with type  $\Gold_5^{\,1}$. The action type for this graph is $(\Gold_5)$.
\medskip

\noindent {\bf (o) Type  $\Gold_5^{\,1}$}:  \\[+4pt]
The smallest example of type  $\Gold_5^{\,1}$ is a graph of order 990 
having automorphism group $\Aut({\rm M}_{12})$ of order 190080, with action type $(\Gold_5,\Gold_5^{\,1})$.  
This is the only example of order up to 10000. 
\medskip

It remains an open question to find the smallest examples of connected finite semisymmetric cubic graphs 
of types $\Gold_2^{\,2}$, $\Gold_2^{\,3}$, $\Gold_3$, $\Gold_4$ and  $\Gold_5$. 
It would also be interesting to know precisely which action types are possible in the semisymmetric case, 
as was done in \cite{CN} for the arc-transitive case.

\section{Existence of an infinite family for each amalgam type}
\label{sec:Existence}

Before continuing, we provide some theoretical background about graph coverings, and
lifting groups of automorphisms along covering projections.
  
Let $\tGa$ and $\Gamma$ be two connected simple graphs, and let $\wp\colon \tGa\to \Gamma$ 
be a graph morphism -- that is, a function from $\V(\tGa)$ to $\V(\Gamma)$  
such that $\wp(u) \sim_\Gamma \wp(v)$ whenever $u\sim_\tGa w$.
 
If the restriction of $\wp$ to the neighbourhood $\tGa(\tv)$ is a bijection between $\tGa(\tv)$ and $\Ga(\wp(\tv))$, 
for every $\tv\in \V(\tGa)$, then we say that $\wp$ is a {\em covering projection}, 
and in that case we call the pre-image $\wp^{-1}(v)$ of a vertex $v$ of $\Ga$ a {\em vertex-fibre}.
For any such $\wp$, the group of all automorphisms of $\tGa$ that preserve every vertex-fibre set-wise is called the {\em group of covering transformations} of $\wp$, and is denoted by $\CT(\wp)$.  
It is easy to see that $\CT(\wp)$ acts faithfully and indeed semiregularly on each vertex-fibre.
 If it acts  transitively (and hence regularly) on each fibre, then the covering projection $\wp$ is said to 
 be {\em regular}.
  
Next, if $\tg\in \Aut(\tGa)$ and $g\in \Aut(\Ga)$ satisfy $\wp(\tv^\tg) = \wp(\tv)^g$ for every $\tv\in \V(\tGa)$, 
then we say that $g$ {\em lifts\/} along $\wp$, and that $\tg$ is a {\em lift\/} of $g$, and also  
that $\tg$ {\em projects\/} along $\wp$ and that $g$ is a {\em projection\/} of $\tg$.
Moreover, if $G\le \Aut(\Ga)$ and every $g\in G$ lifts along $\wp$, then the set of all lifts of elements of $G$ 
is called the {\em lift\/} of $G$, and similarly, if every element of $\tG \le \Aut(\tGa)$ projects, then the set 
of all projections of elements of $\tG$ is called the {\em projection\/} of $\tG$.
  
With this terminology, the group $\CT(\wp)$ is precisely the lift of the trivial subgroup of $\Aut(\Ga)$, 
and  a lift of a projection of some group $\tG\le\Aut(\tGa)$ is the group $\CT(\wp)\tG$.
It is also easy to see that $\tG \le \Aut(\tGa)$ projects along $\wp$ if and only if it normalises the 
covering transformation group $\CT(\wp)$.

\medskip

Next let $\Gamma$ be a graph, with arc-set $A=\A(\Gamma)$, and cycle-set ${\mathcal C}(\Gamma)$. 
Also let $\ZZ A$ be the free $\ZZ$-module over $A$, and let $R$ be the sub-module 
$\langle\, x+x^{-1} : x \in A\,\rangle$ of $\ZZ A$, where $x^{-1}$ denotes the reverse of an arc $x$.
For a cycle $C$ that traverses the arcs $x_0=(v_0,v_1)$, $x_1=(v_1,v_2), \ldots,$
 $x_{n-1}=(v_{n-1},v_0)$ in that order, let $a_C$ be the image of  $x_0 + x_1 + \ldots + x_{n-1} \in \ZZ A$ in the quotient module $\ZZ A/R$.

Then the submodule $\langle\, a_C : C \in {\mathcal C(\Gamma)} \,\rangle$ of $\ZZ A$ is called 
the {\em first homology group} (or sometimes also {\em the integral cycle space}) of $\Gamma$, 
and is denoted by $\H_1(\Gamma;\ZZ)$.
Intuitively we may think of $\H_1(\Gamma;\ZZ)$ as the submodule of $\ZZ A$ generated by all the cycles 
of $\Gamma$, where the reverse $x^{-1}$ of each arc $x$ is identified with the quantity $-x \in \ZZ A$.
 
Since every automorphism of $\Gamma$ takes cycles to cycles, there exists a natural action 
of $\Aut(\Gamma)$ on $\H_1(\Gamma;\ZZ)$, preserving the structure of the $\ZZ$-module. 

In very specific circumstances, this action could potentially be unfaithful, but in that case, 
a non-trivial automorphism of $\Gamma$ would have to preserve each cycle of $\Gamma$ 
while acting on it as a rotation (and not as a reflection). 
This happens, for example, when $\Gamma$ itself is a cycle.

When this pathological situation does not arise, the following theorem can be applied to $\Gamma$.

\begin{theorem}\cite[Theorem 6]{PScov}
\label{thm:PScov}
Let $p$ be an odd prime, let $\Gamma$ be a finite connected graph such that the induced action 
of $\Aut(\Gamma)$ on $\H_1(\Gamma;\ZZ)$ is faithful, and  let $G\le \Aut(\Gamma)$.  
Then there exists a regular covering projection $\wp \colon \tGa \to \Ga$, with $\tGa$ finite, 
such that the maximal group that lifts along $\wp$ is $G$, and the group of covering transformations 
of $\wp$ is a $p$-group.
\end{theorem}

Equipped with the above, we can now prove Theorem~\ref{thm:existence} from the Introduction. 
In fact we prove the following more detailed version of Theorem~\ref{thm:existence}, 
in which `a group of amalgam type' means a group acting edge-transitively with Goldschmidt type 
or Djokovi{\' c}-Miller type.

\begin{theorem}
\label{thm:existence2}
Let $\Gamma$ be a finite cubic edge-transitive graph, and let $G$ be an edge-transitive group 
of automorphisms of $\Gamma$ of amalgam type $T$.
Then for all but finitely many primes $p$, there exists a regular covering projection $\wp \colon \tGa \to \Ga$ 
such that $\tGa$ is finite, $\Aut(\tGa)$ is the lift of $G$ along $\wp$, and the group $\CT(\wp)$ of covering 
transformations of $\wp$ is a $p$-group. In particular, also $\tGa$ has type $T$.
\end{theorem}

\begin{proof}
We begin by proving that $\Aut(\Gamma)$ acts faithfully on $\H_1(\Gamma;\ZZ)$. 
This first part of the proof follows the proof of \cite[Lemma 7]{PScov} almost verbatim.
 
 Assume to the contrary that the action of $\Aut(\Gamma)$ on $\H_1(\Gamma;\ZZ)$ is not faithful. 
Then there exists an automorphism $g$ fixing every element of $\H_1(\Ga;\ZZ)$, and a vertex $v$ of $\Gamma$ 
such that $v^g \not = v$.
Now let $u,w$ and $z$ be the three neighbours of $v$. As $\Gamma$ is finite and contains no vertices of valency $1$, it contains a cycle, and then since $\Gamma$ is edge-transitive, every one of its edges lies on 
a cycle. Moreover, the stabiliser $\Aut(\Gamma)_v$ acts transitively on the neighbourhood 
$\Gamma(v) = \{u,w,z\}$ of $v$, so it contains an element that cyclically permutes $u,w$ and $z$, 
and hence each of the three $2$-paths centred at $v$ lies on a cycle.

Let  $C_1$ be a cycle through the $2$-arc $(u,v,w)$, and $C_2$ be a cycle through the $2$-arc $(u,v,z)$, 
then fix an orientation of $C_1$ and $C_2$ in such a way that $u$ is a predecessor of $v$ in both $C_1$ and $C_2$, and consider $C_1$ and $C_2$ as elements of $\H_1(\Ga;\ZZ)$. By assumption, $g$ preserves $C_1$ and $C_2$, as well as their orientations. In particular, the vertex $v^g$ lies on $C_1$ and on $C_2$, 
and if we let $P_i$ be the path from $v^g$ to $u$ following the cycle $C_i$ in the positive direction 
with respect to the chosen orientation, for $i\in \{1,2\}$, then since $g$ preserves orientation, the vertex 
$u^g$ lies on neither $P_1$ nor $P_2$. Now if we let $C$ be the closed walk obtained by concatenating 
$P_1$ with the reverse of $P_2$, and consider it as an element of $\H_1(\Ga;\ZZ)$, then by assumption, 
$C$ is fixed by $g$, but the vertex $u$ belongs to $C$, while $u^g$ does not, a contradiction. 
This proves that $\Aut(\Gamma)$ acts faithfully on $\H_1(\Gamma;\ZZ)$, as claimed.

The second part of the proof closely follows the proof of \cite[Theorem 9]{PScov}. 

Recall that the order of the edge-stabiliser in the automorphism group of a finite cubic edge-transitive graph
is bounded above by the constant $c=128$. 
Let $n$ be the number of edges of $\Gamma$, and let $p$ be any prime such that $p >nc$.
Then since $\Aut(\Gamma)$ acts faithfully on $\H_1(\Gamma;\ZZ)$, we may use Theorem~\ref{thm:PScov} 
to obtain a regular covering projection $\wp \colon \tGa \to \Ga$ with $\tGa$ finite, such that the maximal 
group that lifts along $\wp$ is $G$, and the group $\tK = \CT(\wp)$ is a $p$-group. 

Let $\tG$ be the lift of $G$ along $\wp$. Then $\tK$ is a normal $p$-subgroup of $\tG$ 
such that $\tG/\tK \cong G$, 
and as $G$ is the maximal group that lifts along $\wp$, the normaliser $N_{\tA}(\tK)$ of $\tK$ 
in the automorphism group $\tA = \Aut(\tGa)$ is $\tG$.
Let $\te$ be an edge of $\tGa$ and let $e=\wp(\te)$.
Then by the definition of the constant~$c$, it follows that $|\tA_{\te}| \le c$. 
Next, because $\tG$ is transitive on the edges of $\tGa$, we have  $\tA=\tA_{\te}\tG$, 
and hence $|\tA\!:\!\tG| = |\tA_{\te}\!:\!\tG_{\te}| \le |\tA_{\te}| \le c < p$. 
Also $|\tG\!:\!\tK|=|G|=n|G_e|\le nc<p$, and it follows that $|\tA\!:\!\tK| = |\tA\!:\!\tG||\tG\!:\!\tK|$ is not divisible by $p$, and so $\tK$ is a Sylow $p$-subgroup of $\tA$. Then by Sylow theory, the number of Sylow $p$-subgroups of $\tA$ is $|\tA\!:\!N_{\tA}(\tK)| = |\tA\!:\!\tG| < p$ and is congruent to $1$ modulo $p$, so must be $1$, 
and thus $\tA=N_{\tA}(\tK) = \tG$. 

Finally, because $\tG$ is the lift of $G$ along a regular covering projection, it has the same amalgam type as $G$,
and this completes the proof.
\end{proof}

Note that it is easy to find  a finite edge-transitive cubic graph admitting an edge-transitive group of type $T$, 
for every Goldschmidt or Djokovi\'{c}-Miller amalgam type $T$.
A lot of them can be realised in $K_{3,3}$, while others are realised in larger graphs, and indeed we have 
shown in the previous section that for every such $T$ there exists a  finite edge-transitive cubic graph 
whose full automorphism group has type $T$. 
Hence Theorem~\ref{thm:existence2} shows that there are infinitely many finite edge-transitive cubic graphs 
of each type, as required.

\section{Asymptotic enumeration}
\label{sec:Enumeration}

In this section we will say that a function $f\colon \NN \to \RR$ is of type $n^{\log n}$ if there exist positive real constants $a$ and $b$ such that $n^{a\log n} \le f(n) \le n^{b\log n}$ holds for all sufficiently large  integers $n$.
 
The following theorem essentially follows from \cite[Theorem~1]{MSP}, which was stated in 
the more general setting of pro-$p$ groups, and for the sake of completeness, we provide a sketch of its proof. 

\begin{theorem}
\label{thm:GroupGrowth}
Let $G$ be a group containing a free subgroup $F$ of rank $r\ge 2$ as a normal subgroup of finite index, 
let $p$ be a prime not dividing the index $|G\!:\!F|$, and let $f\colon \NN \to \RR$ be the function defined by 
$f(n) = |\cN_n|$, 
where $\cN_n = \{N \le F \mid N\norml G \ and \ |F\!:\!N| = p^\alpha \le n \hbox{ for some integer } \alpha\}$.
 Then $f$ is of type $n^{\log n}$.
\end{theorem}

\begin{proof}
First, $f(n)$ is clearly smaller than the number of all $r$-generated $p$-groups of order at most $n$, 
so an upper bound on $f(n)$ of the form $n^{b\log n}$ for some positive real constant $b$ follows 
from a theorem of Lubotzky \cite[Theorem 1]{Lub}.  
Hence we concentrate on obtaining a lower bound for $f(n)$ of the same form.

Let $F_1 \ge F_2 \ge F_3 \ge \ldots$ be the lower $p$-central series for the group $F$, 
given by $F_1 = F$ and $F_{i+1} = [F,F_i]F_i^{\,p}$ for all $i\ge 1$.
Then for each $i$ the factor group $F_i/F_{i+1}$ is an elementary abelian $p$-group, 
so can be viewed as a vector space $V_i$ over $\GF(p)$. 
Next, define $w_i$ and $s_i$ by $\,|F_i\!:\!F_{i+1}| = p^{w_i}\,$ and $\,|F\!:\!F_i| = p^{s_i}\,$ 
(so that $s_1 = 0$ and $s_i = w_1 + \dots + w_{i-1}$ for $i > 1$). Then $\dim(V_i) = w_i$ for all $i$.
Also define $X_i = \{ N\norml G \mid F_{i+1} \le N \le F_i\}$, so that $f(n) \ge |X_i|$ for all $i$. 

Since all the $F_i$  are characteristic subgroups of $F$ and hence normal in $G$, conjugation 
by $G$ induces an action of $G$ on each factor $F_i/F_{i+1}$, with kernel containing $F$, 
because $F_i/F_{i+1}$ is central in $F/F_{i+1}$. 
Accordingly, this gives a linear representation $G/F\to \GL(V_i)$.
Note also that there is a bijective correspondence between the set of invariant subspaces of this linear representation and the set $X_i$ defined in the previous paragraph.

On the other hand, by \cite[Lemma~1]{MSP}, there exists a positive constant $c$, depending solely 
on the group $G/F_i$ and the prime $p$ (which must be coprime to $|G/F|$),
such that every linear representation of $G/F$ on a $w$-dimensional vector-space over $\GF(p)$ 
has at least $p^{cw^2}$ invariant subspaces.
Together with the observations above, this implies that $|X_i| \ge p^{cw_i^2}$ and hence that
$f(n) \ge p^{cw_i^2}$, for all $i$. 

Now define $a = \frac{4c}{9r^2}$ (where $r$ is the rank of the free subgroup $F$), 
and let $n$ be any positive integer such that $n \ge |F\!:\!F_2| = p^{s_2}$,  
and for this $n$, let $j\ge 1$ be such that $p^{s_{j+1}} \le n < p^{s_{j+2}}$.

We will prove that $3r w_j \ge 2 s_{j+2},$ from which it follows that $9r^2 w_j^{\,2} \ge 4 (s_{j+2})^2,$ 
and hence that  \\[+2pt]
$cw_j^{\,2} \ge \frac{4c}{9r^2} (s_{j+2})^2 = a (s_{j+2})^2 > a (\log_p n)^2,$ 
and so $f(n) \ge p^{cw_{j}^{\,2}}  \ge  p^{a (\log_p n)^2} = n^{a \log_p n} = n^{(\frac{a}{\log p}) \log n},$ \\[+2pt]
giving us a lower bound on $f(n)$ of the form we require.
\smallskip
 
To do this, we derive estimates of the values of the parameters $w_j$ and $s_{j+2}$.

It follows from the work of Bryant and Kovacs \cite{BK} that the dimension $w_j$ of
the linear space $F_{j}/F_{j+1}$ is $\sum_{i=1}^{\ j} r_i$, where $r_i$ is the
rank of the $i\,$th section $\gamma_i(F)/\gamma_{i+1}(F)$ of the lower central series of $F$, 
viewed as a $\ZZ$-module.
The rank $r_i$ was determined by Witt \cite{Witt} and equals $\frac{1}{i}\sum_{\,d\,\mid\, i} \mu(d) r^{i/d}$,
 where $\mu$ denotes the arithmetic M\"obius function (see \cite[Lemma 20.6(iii)]{BNV}). 
 Accordingly
$$
w_j= \sum_{i=1}^{\ j}\frac{1}{i}\sum_{d\,\mid\, i} \mu(d) r^{\frac{i}{d}} \quad \
\hbox{ and } \quad s_{j+2} = \sum_{i = 1}^{j+1} w_i.
$$

As was proved in \cite[Lemma 20.7]{BNV} using a Stolz-Ces\`ar analogue of L'H\^ospital's rule,
 the values of $w_j$ and $s_{j+2}$ can be  asymptotically approximated as
$$
 w_j = \frac{r^{\,j+1}}{j(r+1)}(1+o(1)), \quad \hbox{ and } \quad  s_{j+2} = \frac{r^{\,j+3}}{(j+2)(r+1)^2}(1+o(1)), \quad \hbox{ as } j \to \infty.
$$
${}$\\[+2pt]
In particular, as $\frac{2}{3} \le \frac{r}{r+1} < 1$ and $\frac{1}{j+2} <\frac{1}{j}$, 
it follows that for sufficiently large values of $j$, we have 
$\,w_j \ge \frac{2r^{\, j}}{3j}\,$ and $\,s_{j+2} \le  \frac{r^{\,j+1}}{j},\,$ 
which immediately gives us $\,3r w_j \ge 2 s_{j+2},\,$ 
completing the proof.
\end{proof}

The above theorem has the following straightforward consequence.

\begin{corollary}
\label{cor:GroupGrowth}
Let $G$ be a group containing a free subgroup $F$ of rank $r\ge 2$ as a normal subgroup of finite index, 
and let $p$ be a prime not dividing the index $|G\!:\!F|$. 
Also for every positive integer $n$, let $\cN_n$ be the set defined in Theorem~\ref{thm:GroupGrowth}, 
and for subgroups $K$ and $N$ in $\cN_n$, write $K \sim N$ whenever $G/K$ is isomorphic to $G/N$, and
 let $\,\cN_n^*\,$ denote the partition of $\cN_n$  into the equivalence classes of the relation $\sim$.
 Then the function $g \colon \NN \to \RR$ given by $g(n) = |\cN_n^*|$ is of type $n^{\log n}$.
\end{corollary}

\begin{proof}
Since $F$ has finite rank and finite index in $G$, the group $G$ is generated by some finite set $X$, 
say of size $d$. 
Now for subgroups $K$ and $N$ in the same equivalence class of $\sim$,
let $\pi \colon G \to G/K$ be the natural homomorphism,  and let $\varphi \colon G/K \to G/N$ be an isomorphism. Then the (left-to-right) composite homomorphism $\pi\varphi \colon G \to G/N$ has kernel $K$. 
On the other hand, every homomorphism from $G = \langle X \rangle$ to $G/N$ is determined by the images 
of the elements of $X$, and hence there are at most $|G/N|^d$ such homomorphisms, and hence at most 
$|G/N|^d$ possibilities for $K$. Finally, if $h = |G\!:\!F|$, then $|G/N| = |G/F||F/N| = hn$, and it follows that the 
the equivalence class of $N$ under $\sim$ has size at most $(hn)^d$, 
and therefore $g(n) \ge \frac{f(n)}{h^{d} n^{d}}$. Then since $f$ is of type $n^{\log n}$, so is $g$.
\end{proof}

We can now strengthen Theorem~\ref{thm:existence} to the following one. 

\begin{theorem}
\label{thm:enumeration}
Let $\Gamma$ be a connected finite cubic graph, and let $G$ be an edge-transitive group 
of automorphisms of $\Gamma$.
Then there exists a positive real constant $a$ such that for all sufficiently large $n$, 
the number of (non-isomorphic) finite regular coverings of $\tGa$ with $|V(\tGa)| \le n$
whose automorphism group is precisely the lift of $G$ is at least $n^{a \log n}$.
\end{theorem}

\begin{proof}
By Theorem~\ref{thm:existence2}, there exists a regular covering projection $\wp \colon \tGa \to \Ga$ 
such that $\Aut(\tGa)$ is the lift of $G$.
Now let $\wp_\circ \colon T_3 \to \tGa$ be the universal covering projection from the cubic tree $T_3$ to $\tGa$,  and let $A$ be the lift of $\Aut(\tGa)$ along $\wp_\circ$. Then $\Aut(\tGa) \cong A/F$ for some free normal subgroup $F$ of $A$, acting semiregularly on the vertices and on the edges of $T_3$.
Note that the rank of $F$ is equal to the rank of the group $\H_1(\tGa,\ZZ)$, which can be computed 
as $|V(\tGa)| - |E(\tGa)| + 1$.
Moreover $\wp_\circ$ is equivalent 
to the quotient projection $T_3 \to T_3/F$. 
Hence we may assume that $\tGa = T_3/F$ and $\Aut(\tGa) = A/F$.

Next, let $p$ be a prime not dividing $|A\!:\!F|$, and as before (but with $A$ taking the role of $G$), 
for any given positive integer $n$, 
let $\cN_n = \{N \le F \mid N\norml A \ \hbox{and} \ |F\!:\!N| = p^\alpha \le n \hbox{ for some } \alpha\}$.
Then by Corollary~\ref{cor:GroupGrowth}, there exists a positive real constant $a$ 
and a subset $\cN_n^*\subseteq \cN_n$ such that $|\cN_n^*| \ge n^{a\log n}$ for every large enough $n$, 
with the property that for all $K$ and $N$ in $\cN_n^*$, if $A/K \cong A/N$ then  $K=N$.

Now take some $N \in \cN_n^*$ and consider the quotient graph $\tGa_N=T_3/N$. 
Note that since $N\le F$, the quotient projection  $\wp_N \colon T_3 \to \tGa_N$ is a regular covering projection.
Moreover, the group $F/N$ acts semiregularly on the vertices and edges of $\tGa$, 
and so the quotient projection $\wp_{F/N} \colon \tGa_N \to \tGa_N/(F/N)$ is also a covering projection. 
But $\tGa_N/(F/N)$ is isomorphic to $(T_3/N) / (F/N) \cong T_3/F = \tGa$ in a natural way, 
and so in this sense, we may conclude that $\wp_\circ = \wp_{F/N} \circ \wp_N$.
Furthermore,  the group $\Aut(\tGa) = A/F$ lifts along $\wp_{F/N}$ to the group $A/N$.
Then since the edge-stabiliser for a connected finite cubic edge-transitive graph has order at most $128$, 
it follows that $|\Aut(\tGa_N)| \le 128 |E(\tGa_N)| \le 128 |A/N|$. 
In particular, if $F/N$  is a $p$-group for some  prime $p > 128 |(A/N)\!:\!(F/N)| = 128 |A\!:\!F| = 128 |\Aut(\tGa)|$,
then (just as in the proof of Theorem~\ref{thm:existence2}), $F/N$ is a normal Sylow $p$-subgroup 
of $\Aut(\tGa_N)$, and therefore $\Aut(\tGa_N)$ projects via the quotient projection $\tGa_N \to \tGa$, 
so $\Aut(\tGa_N)$ is equal to the lift of $\Aut(\tGa)$.

Thus we have shown that every $N \in \cN_n^*$ yields a regular covering $\tGa_N$ of $\tGa$ 
(and hence also of $\Gamma$) such that $\Aut(\tGa_N)$ is the lift of $\Aut(\tGa)$ 
along $\wp_{F/N} \colon \tGa_N \to \tGa$ 
(and hence also of $G$ along $\wp \circ \wp_{F/N} \colon \tGa_N \to \Gamma$). 
On the other hand, as $\Aut(\tGa_K) = A/K \not \cong A/N = \Aut(\tGa_N)$ for two distinct 
groups $K$ and $N$ in $\cN_n^*$, it follows that $|\cN_n^*|$ is bounded above by the number of regular
coverings $\tGa'$ of $\tGa$ with $|V(\tGa')| \le n |V(\tGa)|$ whose automorphism group is the lift of $G$. 
Finally, because the function $n\mapsto |\cN_n^*|$ is of type $n^{\log n}$, we find that for some 
positive real constant $a$ and for all sufficiently large $n$, the number of these coverings is 
at most $n^{a \log n}$.
\end{proof}

\begin{remark}
As every Goldschmidt or Djokovi\'{c}-Miller amalgam $T$ can be realised in some finite graph, 
and the amalgam type of a lift of a group $G$ of automorphisms is the same as for $G$,
a consequence of Theorem~\ref{thm:enumeration} is that there exists a positive real constant $a$ such that the number of cubic edge-transitive graphs of order at most $n$ with automorphism group of type $T$ is at least $n^{a\log n}$.
In fact, since every cubic edge-transitive graph $\Gamma$ of order $n$ is uniquely determined 
by its type $T$ and the epimorphism from the corresponding universal group $A \le \Aut(T_3)$ of type $T$ 
to $\Aut(\Gamma)$, it follows that the number of cubic edge-transitive graphs of type $T$
with order up to $n$ is bounded above by the number of normal subgroups of index at most $n$ in $A$ times
the maximum number of epimorphisms from $A$ to a finite group of order at most $n$.
By the same theorem of Lubotzky 
as used earlier, the first factor is at most $n^{b\log n}$ for a constant $b$ depending
only on the minimum number $d$ of generators for $A$, while the second factor is at most $n^d$.
Hence the growth rate of the number of connected edge-transitive cubic graphs of any given
Goldschmidt or Djokovi\'{c}-Miller amalgam type $T$ with order up to $n$ is $n^{\log n}$.

\end{remark}

\bigskip

\noindent
\centerline{\sc Acknowledgements}
\\[+5pt] 
The first author is grateful for support from New Zealand's Marsden Fund (under contract UOA2320), 
and to the second author for his hospitality when much of this work was undertaken.  
The second author acknowledges support of the Slovenian Research and Inovation Agency (ARIS), project numbers P1-0294 and J1-4351.
Both authors acknowledge the valuable help of the {\sc Magma} system \cite{Magma} in finding 
and analysing examples (and patterns among those examples), and in formulating and testing conjectures. 
We would also like to thank Gregor Poto\v{c}nik for his help in setting up and maintaining the online version of the census \cite{OnlineCensus}, 
using the `nauty' package \cite{nauty}, and preparing the drawings in this paper.


\end{document}